\documentclass{elsarticle}

\usepackage{amsmath}
\usepackage{amsthm}
\usepackage{amsopn}
\usepackage{amssymb}
\usepackage[all]{xy}
\usepackage{graphicx}

\allowdisplaybreaks[1]

\newcommand{\mr}[1]{\mathrm{#1}}

\newcommand{\abs}[1]{\lvert #1 \rvert}

\newcommand{\bra}[1]{\langle #1 \rangle}

\newcommand{\ZZ}{\mathbb{Z}}

\newcommand{\QQ}{\mathbb{Q}}
\newcommand{\WW}{\mathbb{W}}
\newcommand{\FF}{\mathbb{F}}
\newcommand{\GG}{\mathbb{G}}
\newcommand{\MS}{\mathbb{S}}

 \newtheorem{thm}[equation]{Theorem}
 
 \newtheorem{lem}[equation]{Lemma}
 \newtheorem{prop}[equation]{Proposition}

 \newdefinition{defn}[equation]{Definition}
 \newdefinition{ex}[equation]{Example}
 \newdefinition{exs}[equation]{Examples}
 \newdefinition{rmk}[equation]{Remark}
\newproof{conventions}{Conventions}
\newproof{Ack}{Acknowledgments}

\numberwithin{equation}{section}
\numberwithin{figure}{section}

\DeclareMathOperator{\Aut}{Aut}

\DeclareMathOperator*{\holim}{holim}

\title{The homotopy groups of $S_{E(2)}$ at $p \ge 5$ revisited}
\author[mjb]{Mark Behrens}
\ead{mbehrens@math.mit.edu}

\address[mjb]{Department of Mathematics, MIT, 77 Massachusetts Avenue, Cambridge, MA  02139}

\begin{document}

\begin{abstract}
We present a new technique for analyzing the $v_0$-Bockstein spectral sequence
studied by Shimomura and Yabe. Employing this technique, we derive a
conceptually simpler presentation of the homotopy groups of the $E(2)$-local
sphere at primes $p \ge 5$. We identify and correct some errors in the original
Shimomura-Yabe calculation. We deduce the related $K(2)$-local homotopy groups,
and discuss their manifestation of Gross-Hopkins duality.
\end{abstract}

\maketitle

\section{Introduction}

The chromatic approach to computing the $p$-primary stable homotopy groups of spheres relies on analyzing the chromatic tower:
$$ \cdots \rightarrow S_{E(2)} \rightarrow S_{E(1)} \rightarrow S_{E(0)}. $$
By the Hopkins-Ravenel chromatic convergence theorem \cite{HopkinsRavenel}, the homotopy inverse limit of this tower is the $p$-local sphere spectrum.  The monochromatic layers are the homotopy fibers given by
$$ M_nS \rightarrow S_{E(n)} \rightarrow S_{E(n-1)}. $$
The associated \emph{chromatic spectral sequence} takes the form
$$ \pi_k M_nS \Rightarrow \pi_{k} S_{(p)}. $$

The quest to understand this spectral sequence was begun by Miller, Ravenel, and Wilson \cite{MillerRavenelWilson}, who observed that the monochromatic layers $M_nS$ could be accessed by the Adams-Novikov spectral sequences
\begin{equation}\label{eq:ANSS}
 H^{s,t}(M_0^n) \Rightarrow \pi_{t-s-n}(M_nS)
\end{equation} 
which, for $p \gg n$,  collapse (e.g. for $n = 2$ this spectral sequence collapses for $p \ge 5$).  
The algebraic monochromatic layers $H^{s,t}(M_0^n)$ may furthermore be inductively computed via $v_k$-Bockstein spectral sequences (BSS)
\begin{equation}\label{eq:vkBSS}
 H^{s}(M_{k+1}^{n-k-1}) \otimes \FF_p[v_k]/(v_k^\infty) \Rightarrow H^s(M^{n-k}_k ). 
\end{equation}
The groups $H^*(M_n^0)$, by Morava's change of rings theorem, are isomorphic to the cohomology of the Morava stabilizer algebra.  Miller, Ravenel, and Wilson computed $H^*(M_0^n)$ at all primes for $n \le 1$ and computed $H^0(M_0^2)$ for $ p \ge 3$.

Significant computational progress has been made since \cite{MillerRavenelWilson}, most notably by Shimomura and his collaborators.  A complete computation of $H^*(M^2_0)$ (and hence of $\pi_*S_{E(2)}$) for $p \ge 5$ was achieved by Shimomura and Yabe in \cite{SY}.  Shimomura and Wang computed $\pi_*S_{E(2)}$ at the prime $3$ \cite{ShimomuraWang}, and have computed $H^*(M_0^2)$ at the prime $2$ \cite{ShimomuraWang2}.  These computations are remarkable achievements.

It has been fifteen years since Shimomura and Yabe published their computation of $\pi_*S_{E(2)}$ for primes $p \ge 5$ \cite{SY}.  Since this computation, many researchers have focused their attention on $v_2$-periodic phenomena at ``harder primes'', most notably at the prime $3$, regarding the generic case of $p \ge 5$ as being solved.
Nevertheless, the author has been troubled by the fact that while the image of the $J$-homomorphism ($\pi_* S_{E(1)}$) is familiar to most homotopy theorists, and the Miller-Ravenel-Wilson $\beta$-family ($H^0(M_0^2)$) is well-understood by specialists, the Shimomura-Yabe calculation of $\pi_*S_{E(2)}$ is understood by essentially \emph{nobody} (except the authors of \cite{SY}).  Perhaps even more troubling to the author was that even after careful study, he could not conceptualize the answer in \cite{SY}.  In fact, the author in places could not even parse the answer.  

The difficulties that the author reports above regarding the Shimomura-Yabe calculation (not to mention the Shimomura-Wang computations) might suggest that a complete understanding of the second chromatic layer is of a level of complexity which exceeds the capabilities of most human minds.  However, Shimomura's computation of $H^*(M_1^1)$ (and thus $\pi_*M(p)_{E(2)}$) for $p \ge 5$  \cite{Shimomura} is in fact \emph{very} understandable, and Hopkins-Mahowald-Sadofsky \cite{Sadofsky} and Hovey-Strickland \cite{HoveyStrickland} have even offered compelling schemas to aid in the conceptualization of this computation.  It should not be the case that $\pi_*S_{E(2)}$ is so incomprehensible when the computation of $\pi_*M(p)_{E(2)}$ is so intelligible.

Seeking to shed light on the work of Shimomura-Wang at the prime $3$, Goerss, Henn, Karamanov, Mahowald, and Rezk have constructed and computed with a compact resolution of the $K(2)$-local sphere \cite{GHMR}, \cite{HennKaramanovMahowald}.  Henn has informed the author of a clever technique involving the \emph{projective Morava stabilizer group} that he has developed with Goerss, Karamanov, and Mahowald. When coupled with the resolution, the projective Morava stabilizer group is giving traction in understanding the computation of $\pi_*S_{E(2)}$ at the prime $3$ for these researchers.

The purpose of this paper is to adapt the projective Morava stabilizer group technique to the case of $p  \ge 5$ to analyze the Shimomura-Yabe computation of $\pi_* S_{E(2)}$.  In the process, we correct some errors in the results of \cite{SY} (see Remarks~\ref{rmk:error1}, \ref{rmk:error2}, and \ref{rmk:error3}).  We also propose a different basis than that used by \cite{SY}.  With respect to this basis, $H^*M_0^2$, and consequently $\pi_*S_{E(2)}$ is far easier to understand, and we describe some conceptual graphical representations of the computation inspired by \cite{Sadofsky}.  The author must stress that the errors in \cite{SY} are of a ``bookkeeping'' nature.  The author has found no problems with the actual BSS differentials computed in \cite{SY}.  The computations in this paper are \emph{not} independent of \cite{SY}, as our projective $v_0$-BSS differentials are actually deduced from the $v_0$-BSS differentials of \cite{SY}.

This paper is organized as follows.
In Section~\ref{sec:M02} we review Ravenel's computation of $H^*M_2^0$.  In Section~\ref{sec:M11} we review Shimomura's computation of $H^*M_1^1$ using the $v_1$-BSS.  In Section~\ref{sec:PG} we summarize the projective Morava stabilizer group method introduced by Goerss, Henn, Karamanov, and Mahowald.  This method produces a different $v_0$-BSS for computing $H^*M_0^2$ which we call the \emph{projective $v_0$-BSS}.  In Section~\ref{sec:M20} we show that the differentials in the projective $v_0$-BSS may all be lifted from Shimomura-Yabe's $v_0$-BSS differentials.  We implement this to compute $H^*M_0^2$.  Our computation is therefore not independent of \cite{SY}, but the different basis that the projective $v_0$-BSS presents the answer in makes the computation, and the answer, much easier to understand.  In Section~\ref{sec:SY}, we review the presentation of $H^*M_0^2$ discovered in \cite{SY}, and fix some errors in the process.  We then give a dictionary between our generators and those of \cite{SY}.  In Section~\ref{sec:K(2)S} we review the computation of $\pi_*M(p)_{E(2)}$ and $\pi_*M(p)_{K(2)}$ and give new presentations of $\pi_*S_{E(2)}$ and $\pi_*S_{K(2)}$, using the chromatic spectral sequence.  We explain how these computations are consistent with the chromatic splitting conjecture.  In Section~\ref{sec:GHD} we review the structure of the $K(2)$-local Picard group, and explain how to $p$-adically interpolate the computations of $\pi_*M(p)_{K(2)}$ and $\pi_*S_{K(2)}$.  We explain how Gross-Hopkins duality is visible in $\pi_*M(p)_{K(2)}$.  In Section~\ref{sec:M20new} we give yet another basis for $H^*M_0^2$, which, at the cost of abandoning certain theoretical advantages of the presentation of Section~\ref{sec:M20}, gives an even clearer picture of the additive structure  of $H^*M_0^2$.

\begin{Ack}
It goes without saying that this paper would not have been possible without the previous work of Shimomura and Yabe.  
The author would also like to express his gratitude to Hans-Werner Henn, for explaining the projective Morava stabilizer group method to the author in the first place, to Katsumi Shimomura, for helping the author understand the source of some of the discrepancies found in \cite{SY}, to Tyler Lawson, for pointing out an omission in Lemma~\ref{lem:LHSSS}, and to Paul Goerss and Mark Mahowald, for sharing their $3$-primary knowledge, and helping the author identify a family of errors  in $Y^\infty_1$ and $G^\infty$ in a previous version of this paper.  The author is also very grateful to the time and effort the referee took to carefully read this paper, and provide numerous suggestions and corrections.
The author benefited from the hospitality of the Pacific Institute for Mathematical Sciences and Northwestern University for portions of this work, and was supported by grants from the Sloan Foundation and the NSF.  
\end{Ack}

\begin{conventions}
For the remainder of the paper, $p$ is a prime greater than or equal to $5$.  We define $q$ to be the quantity $2(p-1)$.  We warn the reader that throughout this paper, the cocycle we denote $h_1$ corresponds to what is traditionally called $v_2^{-1}h_1$ (see Section~\ref{sec:M20}).
We will use the notation
$$ x \doteq y $$
to indicate that $x = ay$ for $a \in \FF_p^\times$. 
\end{conventions}

\section{$H^*M_2^0$}\label{sec:M02}

The Morava change of rings theorem gives isomorphisms
$$ H^*(M_2^0) \cong H^*(\GG_2; \pi_* (E_2)/(p, v_1)) \cong H^*(S(2)) \otimes \FF_p[v_2^{\pm 1}] $$
Here $\GG_2$ is the second extended Morava stabilizer group, and $S(2)$ is the second Morava stabilizer algebra.  We refer the reader to \cite{Green} for details.

\begin{thm}[Theorem~3.2 of \cite{Ravenel}]
We have
$$ H^{s,t}(M_2^0) = \FF_p[v_2^{\pm 1}]\{1, h_0, h_1, g_0, g_1, h_0g_1\} \otimes E[\zeta] $$
where the generators have bidegrees $(s,t)$ given as follows.
\begin{align*}
\abs{v_2} & = (0, q(p+1)) \\
\abs{h_0} & = (1,q) \\
\abs{h_1} & = (1,-q) \\
\abs{g_0} & = (2,q) \\
\abs{g_1} & = (2, -q) \\
\abs{\zeta} & = (1,0)
\end{align*}
\end{thm}

\begin{figure}
\includegraphics[width=0.9\textwidth]{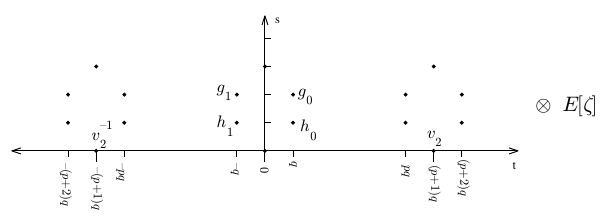}
\caption{$H^* M_2^0$}\label{fig:M02}
\end{figure}

Figure~\ref{fig:M02} displays a chart of this cohomology.

\section{$H^*M_1^1$}\label{sec:M11}

In this section we give a brief account of the structure of the $v_1$-BSS
\begin{equation}\label{eq:v1BSS}
 H^{s}(M_{2}^{0}) \otimes \FF_p[v_1]/(v_1^\infty) \Rightarrow H^s(M^{1}_1 ). 
\end{equation}

We shall use the notation:
\begin{align*}
x_s & := v_2^s x, \quad \text{for} \: x \in H^*M_2^0, \\
G_n & := 
\begin{cases}
v_2^{-p^{n-2}-p^{n-3}- \cdots - 1}g_1,  &  n \ge 1, \\
g_0, & n = 0,
\end{cases} \\
a_n & := \begin{cases}
p^{n-1}(p+1)-1, & n \ge 1, \\
1, & n = 0, 
\end{cases}\\
A_n & := (p^{n-1}+p^{n-2}+\cdots + 1)(p+1).
\end{align*}
Note that $G_1 = g_1$ and $A_0 = 0$.

\begin{thm}[Section~4 of \cite{Shimomura}]
The differentials in the $v_1$-BSS (\ref{eq:v1BSS}) are given as follows:
\begin{align*}
d(1)_{sp^n} & \doteq \begin{cases}
v_1^{a_n} (h_0)_{sp^n-p^{n-1}}, &  n \ge 1, p \not\vert s, \\
v_1(h_1)_s, & n = 0, p \not\vert s,
\end{cases}\\
d(h_0)_{sp^n} & \doteq v_1^{A_n+2} (G_{n+1})_{sp^n}, \quad n \ge 0, s \not\equiv 0, -1 \mod p, \\
d(h_0)_{sp^n-p^{n-2}} & \doteq v_1^{p^{n}-p^{n-2}+A_{n-2}+2}(G_{n-1})_{sp^{n}-p^{n-1}}, \quad n \ge 2, \\
d(h_1)_{sp} & \doteq v_1^{p-1}(g_0)_{sp-1}, \\
d(G_n)_{sp^n} & \doteq v_1^{a_n}(h_0G_{n+1})_{sp^n}, \quad 
n \ge 0, s \not\equiv -1 \mod p.
\end{align*}
The factors involving $\zeta$ satisfy
$$ d(\zeta x) = \zeta d(x). $$
\end{thm}

\begin{figure}
\includegraphics[height=\textwidth, angle=270]{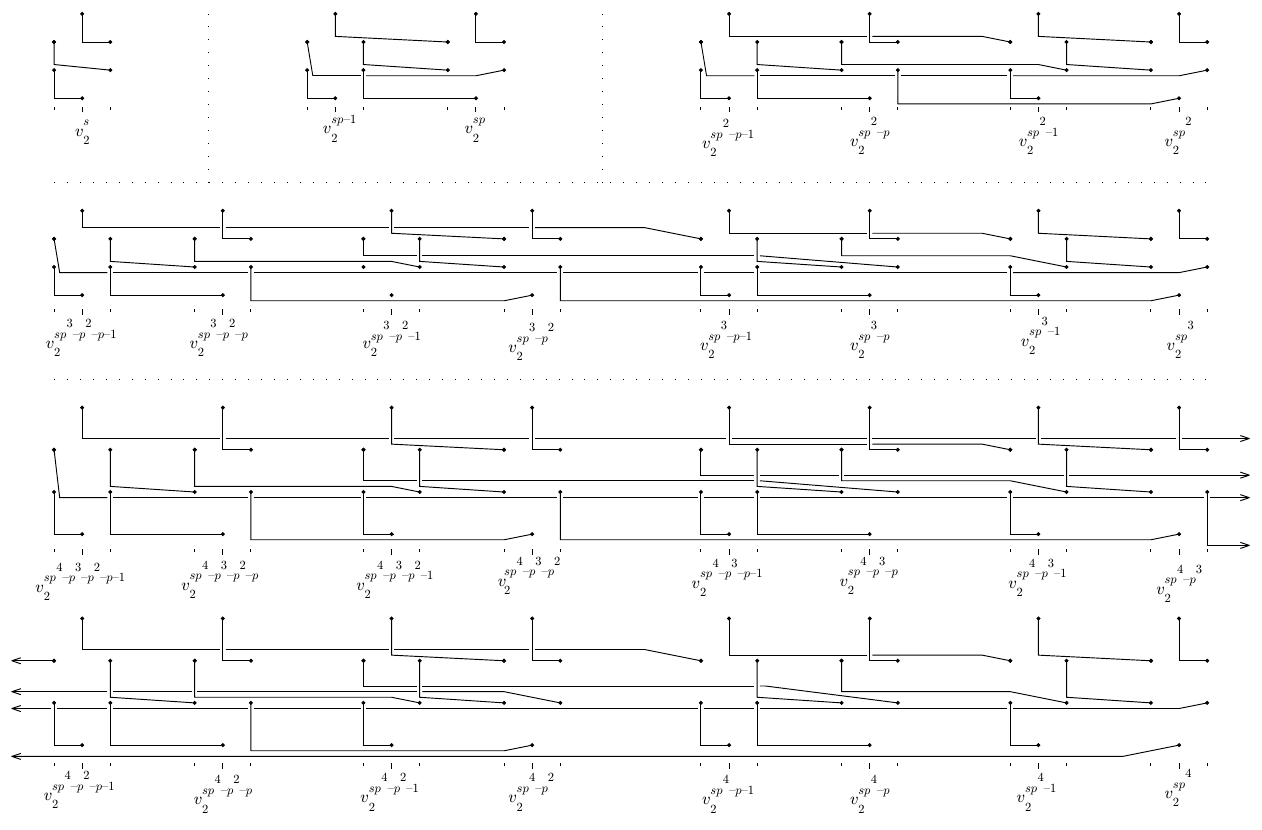}
\caption{$v_1$-BSS in vicinity of $v_2^sp^n$, $0 \le n \le 4$, $s \not\equiv 0, -1 \mod p$, excluding $\zeta$ factor.}\label{fig:v1BSS}
\end{figure}

Figure~\ref{fig:v1BSS} gives a graphical description of these patterns of differentials (excluding the $\zeta$ factors).  In the vicinity of $v_2^{sp^n}$, $s\not\equiv 0, -1 \mod p$, the only elements that are coupled are those of the form
$$ x_{sp^n-\epsilon_{n-1}p^{n-1}- \epsilon_{n-2}p^{n-2} - \cdots - \epsilon_0} $$
for $\epsilon_i \in \{0, 1\}$.

For example, in the vicinity of $v_2^{sp}$, Figure~\ref{fig:v1BSS}  shows the following pattern of differentials.
\begin{center}
\includegraphics{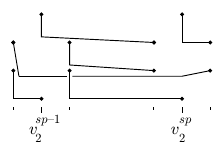}
\end{center}
This depicts the $v_1$-BSS differentials 
\begin{align*}
d(1)_{sp} & \doteq v_1^p (h_0)_{sp-1}, \\
d(1)_{sp-1} & \doteq v_1 (h_1)_{sp-1}, \\
d(h_0)_{sp} & \doteq v_1^{p+3}  (g_1)_{sp-1}, \\
d(h_1)_{sp} & \doteq v_1^{p-1}  (g_0)_{sp-1}, \\
d(g_0)_{sp} & \doteq v_1 (h_0 g_1)_{sp}, \\
d(g_1)_{sp} & \doteq v_1^p (h_0g_1)_{sp-1}. 
\end{align*}
The advantage to using this `hook notation' for the $v_1$-BSS differentials is that the groups $H^* M^1_1$ are easily read off of the diagram.  For example, the hook connecting $(1)_{sp}$ and $(h_0)_{sp-1}$ indicates that there is a $v_1$-torsion summand 
$$ \FF_p[v_1]/(v_1^p)\{ \tfrac{v_2^{sp}}{v_1^p} \} \subset H^0 M^1_1 $$
(generated by $\frac{v_2^{sp}}{v_1^p}$).
Also, the short exact sequence
$$ 0 \rightarrow M_2^0 \xrightarrow{1/v_1} \Sigma^{-q}M^1_1 \xrightarrow{v_1} M^1_1 \rightarrow 0 $$
induces a long exact sequence
$$ \cdots \rightarrow H^s M_2^0 \xrightarrow{1/v_1} H^s M_1^1 \xrightarrow{v_1} H^s M^1_1 \xrightarrow{\delta} H^{s+1}M_2^0 \rightarrow \cdots. $$
The fact that the hook hits $(h_0)_{sp-1}$ indicates that $\delta(\tfrac{v_2^{sp}}{v_1^p}) = (h_0)_{sp-1}$.

The hook patterns of Figure~\ref{fig:v1BSS} can be produced in an inductive fashion.  We explain this inductive procedure below, with a graphical example in the case of $n = 2$.

{\bf Step 1.}  Start with the pattern in the vicinity of $v_2^{sp^{n-1}}$.
\begin{center}
\includegraphics{v1BSSex1.pdf}
\end{center}

{\bf Step 2.} Double the pattern.
\begin{center}
\includegraphics{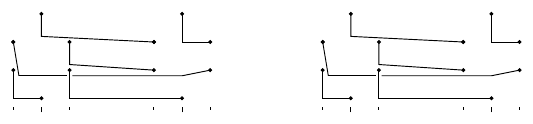}
\end{center}

{\bf Step 3.} Delete the following differentials:
\begin{itemize}
\item the rightmost longest differential on the 0-line,
\item both of the longest differentials on the 1-line,
\item the leftmost longest differential on the 2-line.
\end{itemize}
\begin{center}
\includegraphics{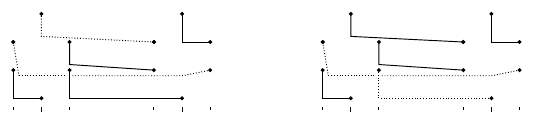}
\end{center}

{\bf Step 4.} Add the following differentials:
\begin{itemize}
\item a differential of length $a_n$ with source $(1)_{sp^n}$,
\item a differential of length $a_n$ with source $(G_n)_{sp^n}$.
\end{itemize}
There are now four elements on the $1$ and $2$ lines left to be connected by differentials.  Couple the closest two, and the farthest two, with differentials.
\begin{center}
\includegraphics{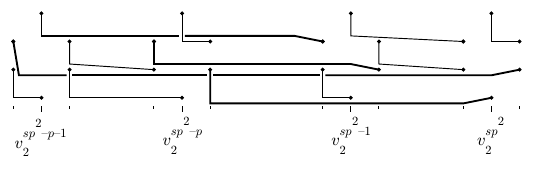}
\end{center}

The cohomology groups $H^*M_1^1$ are easily deduced from the differentials above.  
A complete computation of the groups $H^s(M_1^1)$ first appeared in \cite{Shimomura}.  In that paper, the case of $s=0$ appears as (4.1.5),  and is basically a restatement of the work in \cite{MillerRavenelWilson}.  The case of $s=1$ appears as (4.1.6), and relies on work in \cite{ShimomuraTamura}.  The case of $s > 1$ is covered by Theorem~4.4 of that paper.  Another reference for this result is page 78ff of \cite{HoveyStrickland}, where the translation to the $K(2)$-local setting is given.

The cohomology groups $H^*M^1_1$ are given in Theorem~\ref{thm:M11} below, which uses the notation
$$ x_{s/j} := v_1^{-j} v_2^s x, \quad \text{for} \: x \in H^*M_2^0. $$
However, the reader should be warned, this notation can be misleading, as it is the name of an element in the $E_1$-term of spectral sequence~(\ref{eq:v1BSS}) which detects the corresponding element in $H^* M_1^1$.  For example (c.f.~\cite[p.~190]{Green}) the element $(1)_{p^2/(p^2+1)} \in H^0 M_1^1$ is actually represented by the primitive element
$$ \frac{v_2^{p^2}}{v_1^{p^2+1}} - \frac{v_2^{p^2-p+1}}{v_1^2}- \frac{v_2^{-p} v_3^p}{v_1} \in M_1^1. $$ 

\begin{thm}[\cite{Shimomura}]\label{thm:M11} We have
$$ H^*M_1^1 \cong (X \oplus X_\infty \oplus Y_0 \oplus Y_1 \oplus Y \oplus Y_\infty \oplus G) \otimes E[\zeta]  $$
where:
\begin{align*}
X & := \FF_p\{ 1_{sp^n/j}\}, \quad p \not\vert s, n \ge 0, 1 \le j \le a_n, \\
Y_0 & := \FF_p\{ (h_0)_{sp^n/j}\}, \quad s \not\equiv 0, -1 \mod p, n \ge 0, 1 \le j \le A_n + 2, \\
Y & := \FF_p\{(h_1)_{sp/j}\}, \quad 1 \le j \le p-1, \\
Y_1 & := \FF_p\{(h_0)_{sp^n-p^{n-2}/j}\}, \quad n \ge 2, 1 \le j \le p^{n}-p^{n-2} + A_{n-2}+2, \\
G & := \FF_p\{ (G_n)_{sp^n/j} \}, \quad s \not\equiv -1 \mod p, n \ge 0, 1 \le j \le a_n, \\
X_\infty & := \FF_p\{ 1_{0/j} \}, \quad j \ge 1, \\
Y_\infty & := \FF_p\{ (h_0)_{0/j} \}, \quad j \ge 1. \\
\end{align*}
\end{thm}

\begin{figure}
\includegraphics[height=\textwidth, angle=270]{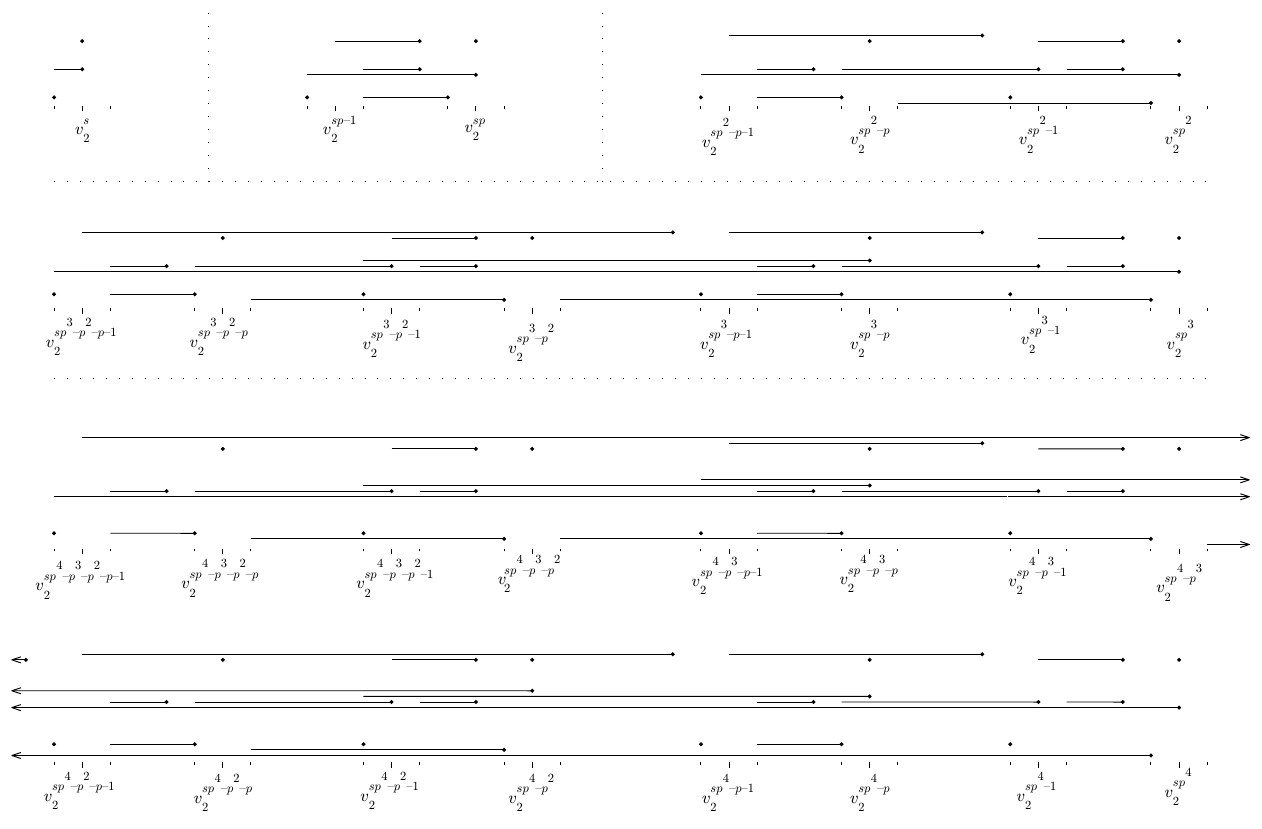}
\caption{$H^* M_1^1$ in the vicinity of $v_2^{sp^n}$, $0 \le n \le 4$, $s \not\equiv 0, -1 \mod p$, excluding the $\zeta$ factor.}\label{fig:M11}
\end{figure}

Figure~\ref{fig:M11} displays pictures of the patterns in this cohomology in the vicinities of $v_2^{sp^n}, s \not\equiv 0, -1 \mod p$ for $0 \le n \le 4$.  The zeta factors are excluded.  In this figure, the patterns are organized according to $v_1$-divisibility.  Thus a family
$$ \FF_p\{ x_{s/j} \}, \quad 1 \le j \le m $$
is represented by:
\begin{center}
\includegraphics{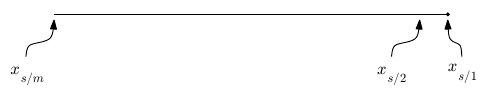}
\end{center}
For example, the pattern in the vicinity of $v_2^{sp}$ depicted in Figure~\ref{fig:M11} is fully labeled below.
\begin{center}
\includegraphics[width=.9\textwidth]{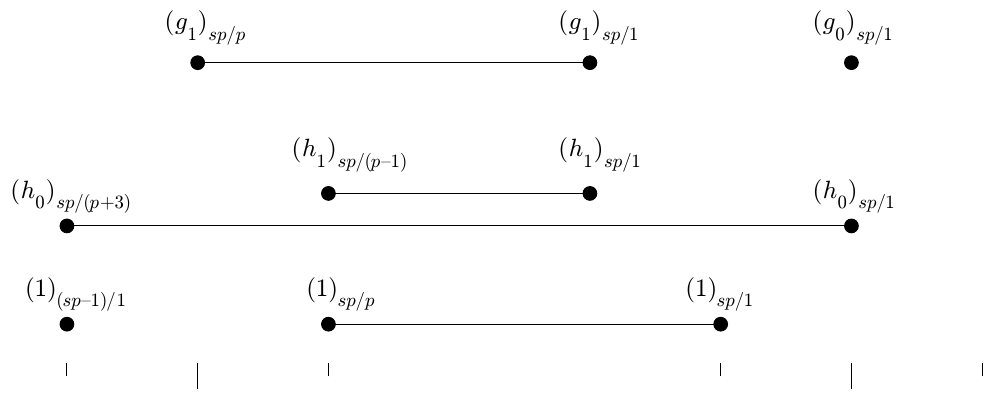}
\end{center}

\section{The projective Morava stabilizer group}\label{sec:PG}

We let $\MS_2$ denote the Morava stabilizer group.  Specifically
$$ \MS_2 := \Aut(H_2) $$
where $H_2$ is the Honda height $2$ formal group over $\FF_{p^2}$.  The action of $\MS_2$ on 
$$ (E_2)_* = W(\FF_{p^2}) [[u_1]] [u^{\pm 1}] $$
extends to an action of the extended Morava stabilizer group
$$ \GG_2 := \MS_2 \rtimes Gal(\FF_{p^2}/\FF_p). $$
Defining 
\begin{align*}
v_1 & := u^{p-1} u_1, \\
v_2 & := u^{p^2-1},
\end{align*}
the Morava change of rings theorem gives isomorphisms:
\begin{align*}
H^*M_2^0 & \cong H^*(\GG_2; (E_2)_*/(p, v_1)), \\
H^*M_1^1 & \cong H^*(\GG_2; (E_2)_*/(p, v_1^\infty)), \\
H^*M_0^2 & \cong H^*(\GG_2; (E_2)_*/(p^\infty, v_1^\infty)).
\end{align*}
We henceforth will use the notation:
\begin{align*}
M_2^0(E_2) & := (E_2)_*/(p, v_1), \\
M_1^1(E_2) & := (E_2)_*/(p, v_1^\infty), \\
M_0^2(E_2) & := (E_2)_*/(p^\infty, v_1^\infty).
\end{align*}

Define the projective (extended) Morava stabilizer group $P\GG_2$ to be the quotient of $\GG_2$ by the center of $\MS_2$.
$$ 1 \rightarrow \ZZ_p^\times \rightarrow \GG_2 \rightarrow P\GG_2 \rightarrow 1. $$

Consider the Lyndon-Hochschild-Serre spectral sequence (LHSSS)
\begin{equation}\label{eq:LHSSS}
H^{s_1}(P\GG_2 ; H^{s_2,t}(\ZZ_p^\times ; M_0^2(E_2)) \Rightarrow H^{s_1+s_2,t}(\GG_2; M_0^2(E_2)).
\end{equation}

The following lemma allow us to analyze (\ref{eq:LHSSS}).

\begin{lem}\label{lem:Zpx}
We have
$$ H^{s,t}(\ZZ_p^\times; M_0^2(E_2)) \cong 
\begin{cases}
[(E_2)_*/(v_1^\infty)]_{t} \otimes \ZZ/p^k, & t = p^{k-1}t'q, p \not\vert t', s = 0, \\
[(E_2)_*/(v_1^\infty)]_{0} \otimes \ZZ/p^\infty, & t = 0, s \in \{ 0, 1\}, \\
0, & \text{otherwise.}
\end{cases}
$$
\end{lem}

\begin{proof}
The subgroup $\ZZ_p^\times \subset \GG_2$ acts on $(E_2)_*$  by the formula
\begin{equation}\label{eq:action} 
[a] \cdot x = a^m x,  \quad a \in \ZZ_p^\times, x \in M_0^2(E_2)_{2m}. 
\end{equation}
The computation is therefore more or less identical to the computation of $H^*M_0^1$.
\end{proof}  

For
$$ \frac{x}{v_1^j} \in [(E_2)_*/(v_1^\infty)]_{t} $$
with $t = p^{k-1}t'q$, we have corresponding elements 
$$ \frac{x}{v_1^j p^k} \in H^{0,t}(\ZZ_p^\times; M_0^2(E_2)). $$
For $x/v_1^j$ in $[(E_2)_*/(v_1^\infty)]_0$ we have elements
\begin{align*}
\frac{x}{v_1^j p^k} & \in H^{0,0}(\ZZ_p^\times; M_0^2(E_2)), \\
\frac{\zeta x}{v_1^j p^k} & \in H^{1,0}(\ZZ_p^\times; M_0^2(E_2)), \\
\end{align*}
for $k \ge 1$.

For dimensional reasons, we deduce the following lemma.

\begin{lem}\label{lem:LHSSS}
For $t \ne 0$, the LHSSS (\ref{eq:LHSSS}) collapses.  In particular, the edge homomorphism (inflation) given by the composite
$$ H^{*,t}(P\GG_2; M_0^2(E_2)^{\ZZ_p^\times}) \rightarrow H^{*,t}(\GG_2; M_0^2(E_2)^{\ZZ_p^\times}) \rightarrow H^{*,t}(\GG_2; M_0^2(E_2)) $$
is an isomorphism for $t \ne 0$.
\end{lem}

\begin{rmk}
Note that the LHSSS (\ref{eq:LHSSS}) also collapses for $t = 0$, though not for dimensional reasons.  See the discussion before Theorem~\ref{thm:M20}.
\end{rmk}

The $p$-adic filtration on $M_0^2(E_2)$ induces a \emph{projective} $v_0$-BSS
\begin{equation}\label{eq:Pv0BSS}
H^{s,t}(P\GG_2; M_1^1(E_2)^{\ZZ_p^\times}) \otimes \FF_p[v_0]/(v_0^{k(t)}) \Rightarrow H^{s,t}(P\GG_2; M_0^2(E_2)^{\ZZ_p^\times})
\end{equation}
where
$$ k(t) := 
\begin{cases}
\nu_p(t)+1, & q \vert t, \\
0, & q \not\vert t.
\end{cases}
$$

The $E_2$-term of (\ref{eq:Pv0BSS}) is easy to understand, as we will now demonstrate. Let $\GG_2^1$ denote the kernel of the reduced norm, given by the composite
$$ \GG_2 \xrightarrow{N} \ZZ_p^\times \rightarrow \ZZ_p^\times/\FF_p^\times \cong \ZZ_p. $$

\begin{lem}
The composite
$$ H^*(P\GG_2 ; M_1^1(E_2)^{\ZZ_p^\times}) \rightarrow H^*(\GG_2 ; M_1^1(E_2)) \rightarrow H^*(\GG_2^1 ; M_1^1(E_2)) $$
is an isomorphism.
\end{lem}

\begin{proof}
Observe there is an isomorphism
$$ P\GG_2 = \GG_2/\ZZ_p^\times \cong \GG^1_2/ (\ZZ_p^\times \cap \GG_2^1) = \GG_2^1/\FF_p^\times. $$
Since $\abs{\FF_p^\times}$ is coprime to $p$, the LHSSS
$$ H^*(P\GG_2; H^*(\FF_p^\times; M_1^1(E_2))) \Rightarrow H^*(\GG_2^1; M_1^1(E_2)) $$
collapses.  Therefore the edge homomorphism gives an isomorphism
$$ H^*(P\GG_2;M_1^1(E_2)^{\FF_p^\times}) \cong H^*(\GG_2^1; M_1^1(E_2)). $$
However, it is immediate from (\ref{eq:action}) that the natural inclusion gives an isomorphism
$$ M_1^1(E_2)^{\ZZ_p^\times} \xrightarrow{\cong} M_1^1(E_2)^{\FF_p^\times}. $$
\end{proof}

The LHSSS
$$ H^*(\ZZ_p; H^*(\GG_2^1; M_1^1(E_2))) \Rightarrow H^*(\GG_2; M_1^1(E_2)) $$
collapses to give an isomorphism
$$
H^*(\GG_2; M_1^1(E_2)) \cong H^*(\GG_2^1; M_1^1(E_2)) \otimes E[\zeta].
$$
The map
$$ H^*(\GG_2; M_1^1(E_2)) \rightarrow H^*(\GG_2^1; M_1^1(E_2)) $$
is the quotient of $H^*(\GG_2; M_1^1(E_2))$ by the zeta factor (see Theorem~\ref{thm:M11}).
We therefore have proven the following lemma.

\begin{lem}\label{lem:PM11}
We have (in the notation of Theorem~\ref{thm:M11}):
$$ H^*(P\GG_2; M_1^1(E_2)^{\ZZ_p^\times}) = X \oplus X_\infty \oplus Y_0 \oplus Y_1 \oplus Y \oplus Y_\infty \oplus G.  $$
\end{lem}

\section{$H^*M_0^2$}\label{sec:M20}

In this section we compute the projective $v_0$-BSS (\ref{eq:Pv0BSS}).  We will deduce our differentials from the differentials of \cite{SY} using the following maps of $v_0$-BSS's.
$$
\xymatrix{
H^{s,t}(P\GG_2; M_1^1(E_2)^{\ZZ_p^\times}) \otimes \FF_p[v_0]/(v_0^{k(t)}) \ar@{=>}[r] \ar[d] &
H^{s,t}(P\GG_2; M_0^2(E_2)^{\ZZ_p^\times}) \ar[d]
\\
H^{s,t}(\GG_2; M_1^1(E_2)) \otimes \FF_p[v_0]/(v_0^{\infty}) \ar@{=>}[r] \ar[d] &
H^{s,t}(\GG_2; M_0^2(E_2)) \ar[d]
\\
H^{s,t}(\GG^1_2; M_1^1(E_2)) \otimes \FF_p[v_0]/(v_0^{\infty}) \ar@{=>}[r] &
H^{s,t}(\GG^1_2; M_0^2(E_2)) 
}
$$
The results of Section~\ref{sec:PG} imply that the composite of these maps on $E_1$-terms is isomorphic to the inclusion
$$ H^{s,t}(\GG^1_2; M_1^1(E_2)) \otimes \FF_p[v_0]/(v_0^{k(t)}) \hookrightarrow H^{s,t}(\GG^1_2; M_1^1(E_2)) \otimes \FF_p[v_0]/(v_0^{\infty}).
$$
The differentials in the middle spectral sequence were computed by \cite{SY}. They therefore map down to differentials in the bottom spectral sequence, and then may be lifted to the top spectral sequence by injectivity.  \emph{In summary: we can regard the $v_0$-BSS differentials of \cite{SY} to be differentials in the projective $v_0$-BSS after we kill all of the terms involving $\zeta$.}

The differentials in the projective $v_0$-BSS (\ref{eq:Pv0BSS}) are given in the theorem below.  Following \cite{SY}, we only list the \emph{leading terms}, which are taken to be the terms of the form $x/v_1^j$ for $j$ maximal.  We will explain why this method suffices in Remark~\ref{rmk:leading}.

\begin{ex}
In Lemma~5.1 of \cite{SY}, it is stated that the connecting homomorphism $\delta: H^0 M_0^2 \rightarrow H^1 M_1^1$ is given on a class $x_{2}/p v_1^{2p} \in M_0^2$ (where $[x_{2}/v_1^{2p}]$ represents $1_{p^{2}/2p} \in H^0M_1^1$) by
$$ \delta(x_{2}/p v_1^{2p}) = -2p y_{p^2}/v_1^{2p+1} - p x_2 \zeta /v_1^{2p} + y_{p^2-1}/v_1^p + v_2^{p^2-p-1}V/v_1^{p-2} + \cdots. $$
Here $[y_s/v_1^j] = (h_0)_{s/j} \in H^1 M_1^1$ and $[v_2^s V/v_1^j] = (h_1)_{s/j} \in H^1 M_1^1$. 
The first two terms are zero, as they have coefficients which are zero mod $p$, but the $\zeta$ term would be ignored anyways for the purposes of the projective $v_0$-BSS.  The leading term is therefore $y_{p^2-1}/v_1^p$, and this corresponds to the projective $v_0$-BSS differential:
$$ d(1_{p^2/2p}) = v_0 (h_0)_{(p^2-1)/p}. $$
\end{ex}

We lift the $v_0$-BSS differentials of \cite{SY} to projective $v_0$-BSS differentials in the following sequence of lemmas.

\begin{lem}
For $p \not\vert s$, $n \ge 0$, $1 \le j \le a_n$, we have:
$$
d(1_{sp^n/j}) \doteq 
\begin{cases}
v_0 (h_0)_{s/2}, & n = 0, j = 1, s \equiv 1 \mod p, \\
v_0(h_1)_{sp/p-1}+ \cdots , & n = 1, j = p, \\
v_0^k (h_0)_{sp^n-p^{n-k-1}/j-a_{n-k}} + \cdots, & n \ge 2, p^k \vert j, a_{n-k} < j \le a_{n-k+1}, \\
0, & \text{in all other cases}.
\end{cases}
$$
We also have
$$ d(1_{0/j}) = 0, \quad j \ge 1.  $$
\end{lem}

\begin{proof}
This follows from Lemma~5.1 of \cite{SY}.  The last assertion is Proposition~6.9(ii) of \cite{MillerRavenelWilson}.
\end{proof}

\begin{lem}
For $1 \le j \le p-1$ we have
$$ d((h_1)_{sp/j}) = 0. $$
\end{lem}

\begin{proof}
This follows from Lemma~7.2 of \cite{SY}.
\end{proof}

\begin{lem}
Let $s \not\equiv 0, -1 \mod p$ and $n \ge 1$.  For $1 \le k \le n$, $A_{n-k}+2 < j \le A_{n-k+1} + 2$, and $p^k \vert j-1$, we have:
$$
d((h_0)_{sp^n/j}) \doteq v_0^k G_{n-k+1/j-A_{n-k}-2} + \cdots.
$$
We have $d(h_0)_{sp^n/j} = 0$ in all other cases.  We also have
$$ d(h_0)_{0/j} = 0,\quad j \ge 1. $$
\end{lem}

\begin{proof}
This follows from Propositions~7.3 and 7.5 of \cite{SY}.  The last assertion follows from the fact that these elements are actually the targets of (non-projective) $v_0$-BSS differentials in Proposition~6.9(ii) of \cite{MillerRavenelWilson}.
\end{proof}

\begin{lem}
Let $n \ge 2$. For $1 \le k \le n-2$, $p^n-p^{n-2}+A_{n-k-2}+2 < j \le p^{n}-p^{n-2} +  A_{n-k-1}+2$, and $p^k \vert j+a_{n-1}$,  we have
$$ d((h_0)_{sp^n-p^{n-2}/j}) \doteq v_0^k(G_{n-k-1})_{sp^n-p^{n-1}/j-p^n+p^{n-2}-A_{n-k-2}-2} + \cdots. $$
We also have
$$ d((h_0)_{sp^n-p^{n-2}/p^{p^n-p^{n-2}+1}}) \doteq v_0^{n-1} (G_0)_{sp^n-p^{n-1}/1}. $$
In all other cases $d((h_0)_{sp^n-p^{n-2}/j}) = 0$.
\end{lem}

\begin{proof}
This follows from Proposition~7.6  of \cite{SY} in the case of $n = 2$, and Proposition~7.8 of \cite{SY} in the case of $n > 2$.  The condition $j > p^n-p^{n-2}+A_{n-k-2}+2$ is not present in Proposition~7.8 of \cite{SY}, but it is necessary because otherwise the target of the differential is not present.
\end{proof}

\begin{figure}
\includegraphics[height=\textwidth, angle=270]{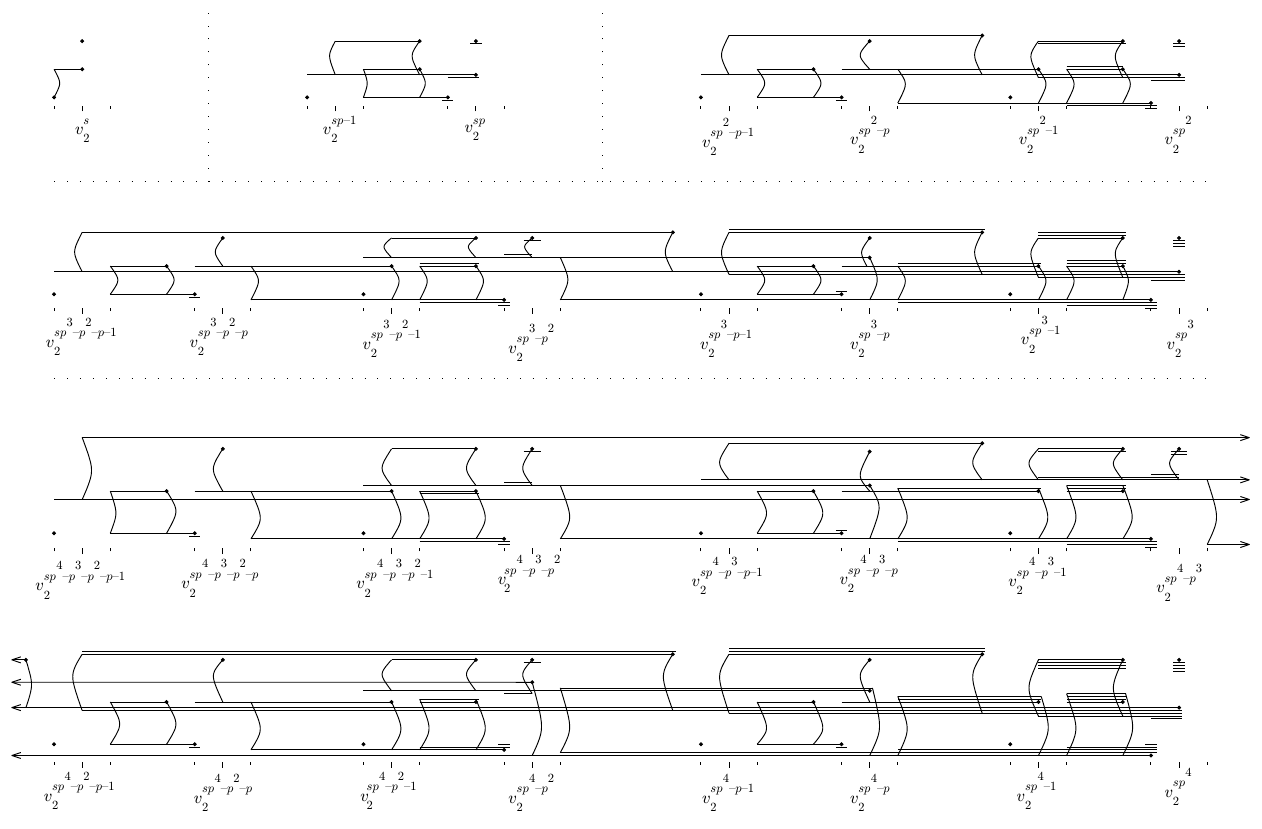}
\caption{$v_0$-BSS in the vicinity of $v_2^{sp^n}$, $0 \le n \le 4$, $s \not\equiv 0, -1 \mod p$.}\label{fig:v0BSS}
\end{figure}

\begin{figure}
\begin{center}
\includegraphics[height=.25\textwidth, angle=270]{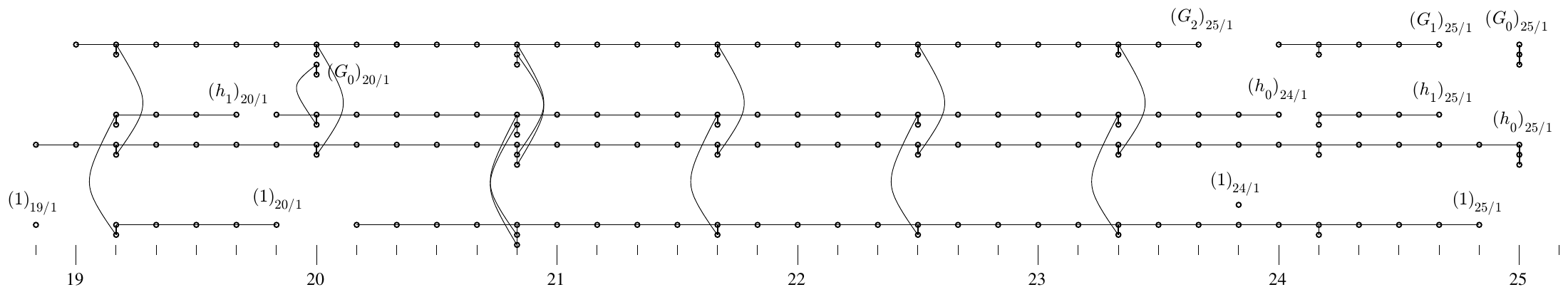}
\hspace{1in}
\includegraphics[height=.25\textwidth, angle=270]{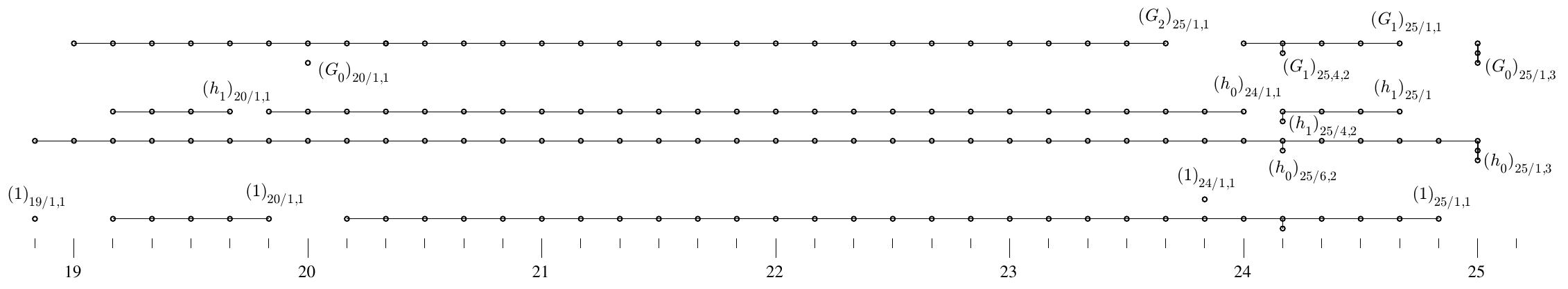}
\end{center}
\caption{Explicit patterns in the case $p = 5$ in the vicinity of $v_2^{25}$: the projective $v_0$-BSS (left) and $H^*M_0^2$ (right).}\label{fig:v0BSSex2-3}
\end{figure}

These theorems account for all of the possible differentials in the projective $v_0$-BSS.
Figure~\ref{fig:v0BSS} displays the patterns of differentials in the projective $v_0$-BSS in the vicinity of $v_2^{sp^n}$, $s \not\equiv 0, -1 \mod p$, for $n \le 4$.  The notation in Figure~\ref{fig:v0BSS} is interpreted as follows.  Given a pair of $k$-fold lines and a region bookended on either side with curved lines as below:
\begin{center}
\includegraphics{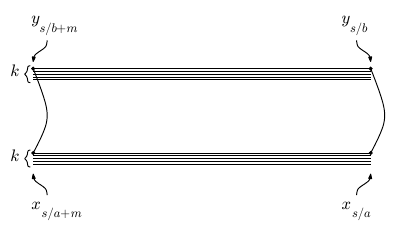}
\end{center}
one has $E_2$-term elements
\begin{align*}
v_0^{-i}x_{s/a+j}, & \quad \text{for} \: 0 \le j \le m, 1 \le i \le \nu_p(\abs{x_{s/a+j}})+1, \\
v_0^{-i}y_{s/b+j}, & \quad \text{for} \: 0 \le j \le m, 1 \le i \le \nu_p(\abs{y_{s/b+j}})+1, \\
\end{align*}
and differentials
$$
d(v_0^{-i}x_{s/a+j}) \doteq v_0^{-i+k} y_{s/b+j} + \cdots,  \quad \text{if} \: \nu_p{\abs{x_{s/a+j}}} \ge k.
$$
Figure~\ref{fig:v0BSSex2-3} shows an explicit example of some of these patterns of differentials in the case where $p = 5$ in the vicinity of $v_2^{25}$.

\begin{rmk}\label{rmk:leading}
The reason it suffices to consider leading terms in the projective $v_0$-BSS differentials is that the differentials are in ``echelon form''.  Firstly, observe that there is  an ordering of the basis of $H^*(P\GG_2; M_1^1(E_2)^{\ZZ_p^\times})$ of Lemma~\ref{lem:PM11} by $v_1$-valuation.  Inspection of the patterns in Figure~\ref{fig:M11} reveal that there are no two basis elements in the same bidegree with identical $v_1$-valuation.  Saying that the projective $v_0$-BSS differentials are in \emph{echelon form} with respect to this ordered basis is equivalent to the assertion that for each $k$, and each pair of elements 
$$ x_{i/j}, x'_{i'/j'} \in H^{s,t}(P\GG_2; M_1^1(E_2)^{\ZZ_p^\times}) $$
with $j < j'$, and with projective $v_0$-BSS differentials
\begin{align*}
d_k(x_{i/j}) & = v_0^k y_{m/l} + \cdots, \\
d_k(x'_{i'/j'}) & = v_0^k y'_{m'/l'} + \cdots,
\end{align*}
we have $l < l'$.  This condition is easily verified to be satisfied by inspecting the patterns in Figure~\ref{fig:v0BSS}.
\end{rmk}

These differentials result in a complete computation of $H^{s,t}(P\GG_2; M_0^2(E_2)^{\ZZ_p^\times})$.  This gives a computation of $H^{s,t}M_0^2$ \emph{except at $t = 0$}.  Using the norm map, one can show that the LHSSS (\ref{eq:LHSSS}) collapses, so that Lemma~\ref{lem:Zpx}  implies that we have
$$ H^{*,0} M_0^2\cong H^{*,0}(P\GG_2; M_0^2(E_2)^{\ZZ_p^\times}) \otimes E[\zeta]. $$ 
In this case the $P\GG_2$ approach offers no advantages over the more traditional $v_0$-BSS:
\begin{equation}\label{eq:t0v0BSS}
H^{*,0}M_1^1 \otimes \FF_p[v_0]/(v_0^\infty) \Rightarrow H^{*,0} M_0^2.
\end{equation}
Moreover Lemma~8.10 of \cite{MillerRavenelWilson}, Corollary~9.9 of \cite{SY}, and Lemma~4.5 of \cite{SY2} imply that there are no non-trivial differentials in (\ref{eq:t0v0BSS}).

We will use the notation
$$ x_{s/j,k} := \frac{v_2^s x}{v_1^j p^k}. $$
Such an element will always have order $p^k$.
The resulting computation of $H^*M_0^2$ is given below.

\begin{thm}\label{thm:M20}
We have 
$$ H^*M_0^2 \cong X^\infty \oplus Y_0^\infty \oplus Y^\infty \oplus Y_1^\infty \oplus G^\infty \oplus X^\infty_\infty \oplus Y_{0,\infty}^\infty \oplus \zeta Y_{0,\infty}^\infty \oplus G_\infty^\infty \oplus \zeta G_\infty^\infty $$
where the summands are spanned by the following elements:
\begin{align*}
X^\infty & := \bra{1_{sp^n/j,k}}, \quad p \not| s, n \ge 0, 1 \le k \le n+1, 1 \le j \le a_{n-k+1}, p^{k-1}|j, \\
X_\infty^\infty & := \bra{1_{0/j,k}}, \quad k \ge 1, j \ge 1, p^{k-1}|j, \\
Y_0^\infty & := \bra{(h_0)_{sp^n/j,k}}, \quad p \not|s, n \ge 0, 1 \le k \le n+1, 1 \le j \le A_{n-k+1}+2, p^{k-1}|j-1, \\
Y_{0,\infty}^\infty  & := \bra{(h_0)_{0/j,k}}, \quad k \ge 1, j \ge 1, p^{k-1}|j-1, \\
\zeta Y_{0,\infty}^\infty & := \bra{\zeta(h_0)_{0/1,k}}, \quad k \ge 1, \\
Y^\infty & :=  \bra{(h_1)_{sp/j,k}}, \quad k = 1, 1 \le j \le p-1, \text{and if} \: p|s, k = 2, j = p-1, \\
Y_1^\infty & := \bra{(h_0)_{sp^n-p^{n-2}/j,k}}, \quad \text{writing} \: s = p^i s', p \not| s', \text{we have:}\\
& \quad  1 \le j \le p^n-p^{n-2}, p^{k-1}|j+a_{n-1}, \: \text{for} \: 1 \le k \le \min(i+1, n+1); \\
& \quad  p^n-p^{n-2} <  j \le p^n-p^{n-2}+A_{n-k-1}+2, p^{k-1} | j+a_{n-1}, \: \text{for} \: 1 \le k \le n-1, \\
G^\infty & := \bra{(G_n)_{sp^n/j,k}}, \quad n \ge 0, 1 \le j \le a_n, \: \text{writing} \: s = p^i t, p \not| t, \: \text{we have}: \\
& \quad \begin{cases}
t \not\equiv -1 \mod p: & i\ge 0, 
\begin{cases}
n = 0: & 1 \le k \le i+1, \\
n \ge 1: & 1 \le k \le \min(n+1, i+1), \\
& p^{k-1}|j+A_{n-1}+1,
\end{cases}
\\
t \equiv -1 \mod p: & i \ge 1,
\begin{cases}
n = 0: & 1 \le k \le i, \\
n \ge 1: & 1 \le k \le \min(n+1, i), \\  
& p^{k-1}|j+A_{n-1}+1,
\end{cases}
\end{cases} \\
G_\infty^\infty & := \bra{(G_n)_{0/j,k}}, \quad n \ge 0, 1 \le j \le a_n, 
\begin{cases} 
n = 0: & k \ge 1, \\
n > 0: & 1 \le k \le n+1, 1 \le j \le a_n, \\
& p^{k-1}|j+A_{n-1}+1,
\end{cases}\\
\zeta G_\infty^\infty & := \bra{\zeta(G_0)_{0/1,k}},\quad k \ge 1. \\
\end{align*}
\end{thm}

\begin{rmk}
Take note that in the theorem above, we have elected to enumerate \emph{all} of the values of $k$ so that the elements $x_{s/j,k}$ exist, not just the maximal values of $k$, which would give a basis.  The author finds that this makes the conditions on the different indices somewhat easier to digest.  The presentation above  does give a basis for the associated graded of $H^*M_0^2$ with respect to the $p$-adic filtration.
\end{rmk}

\begin{figure}
\includegraphics[height=\textwidth, angle=270]{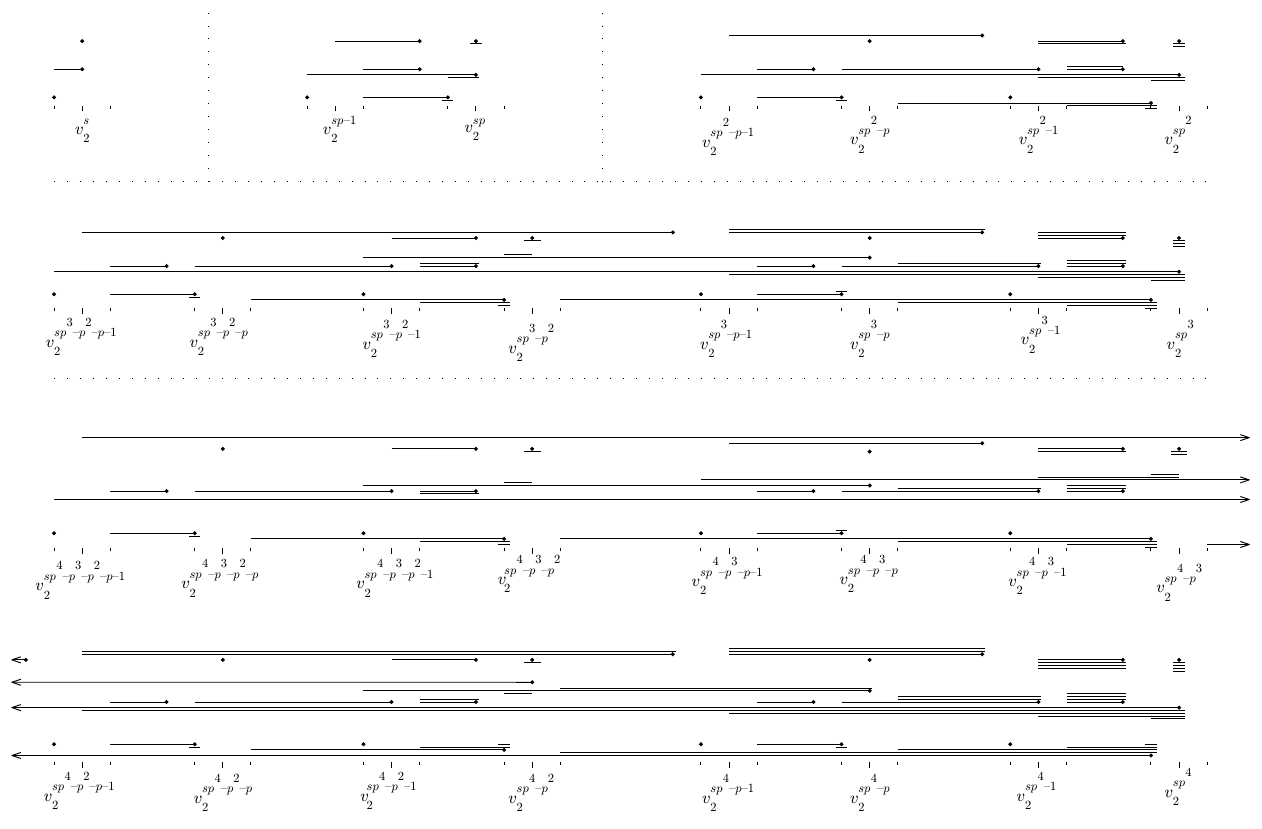}
\caption{$H^* M_0^2$ in the vicinity of $v_2^{sp^n}$, $0 \le n \le 4$, $s \not\equiv 0, -1 \mod p$}\label{fig:M20}
\end{figure}

Figure~\ref{fig:M20} displays the resulting cohomology $H^*M_0^2$ in the vicinities of $v_2^{sp^n}$, $s \not\equiv 0, -1 \mod p$, $n \le 4$.  In this figure, a $k$-fold line segment
\begin{center}
\includegraphics{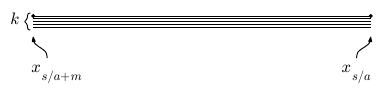}
\end{center}
is spanned by
$$ \bra{x_{s/j,\ell}}, \quad \text{for} \: a \le j \le a+m, 1 \le \ell \le \min(\nu_p(\abs{x_{s/j}})+1, k). $$
Figure~\ref{fig:v0BSSex2-3} shows examples of these patterns in the case where $p=5$ in the vicinity of $v_2^{25}$.

\section{Dictionary with Shimomura-Yabe}\label{sec:SY}

The computation of Shimomura-Yabe uses the $v_0$-BSS
\begin{equation}\label{eq:v0BSS}
H^{s,t}(M_1^1) \otimes \FF_p[v_0]/(v_0^{\infty}) \Rightarrow
H^{s,t}(M_0^2)
\end{equation}
where $H^*(M_1^1)$ is computed as in Theorem~\ref{thm:M11}.
Part of the reason that the computation of $H^*M_0^2$ is so complicated when using this spectral sequence is that the families of Theorem~\ref{thm:M20} get split between families involving $\zeta$ and not involving $\zeta$.  We recall the result of \cite{SY}, with some corrections to their families.  In order to not confuse their generators coming from $H^*(\GG_2; M^1_1(E_2))$ with ours coming from $H^*(P\GG_2; M_1^1(E_2)^{\ZZ_p^\times})$, we will write the Shimomura-Yabe generators, as well as the Shimomura-Yabe families, in non-italic typeface.  We continue to use our $x_{s/j,k}$ notation from Section~\ref{sec:M20}.  We also continue our convention that $\abs{\mr{h}_1} = -q$.

Below we reproduce the main result of \cite{SY}.  Our reason for reproducing the whole answer is that the author could not fully parse the conditions as printed in \cite{SY}.  Also, the author discovered some errors in the paper: the answer below includes the author's corrections.

\begin{thm}[Theorem~2.3 of \cite{SY}]\label{thm:SY}
The cohomology $H^*M_0^2 $ is isomorphic to
\begin{gather*}(\mr{X}_\infty^\infty \oplus \mr{Y}^\infty_{\infty,C} \oplus \mr{G}_0^\infty) \otimes E[\zeta] \oplus \mr{X}^\infty \oplus \mr{X}\zeta^\infty_C  \oplus \mr{Y}^\infty_{0,C} \oplus \mr{Y}^\infty_{1,C} \oplus \mr{Y}^\infty_C \oplus \\
\mr{G}^\infty_C \oplus (\mr{Y}^{\infty,G}_{0,C} \oplus  \mr{Y}^{\infty, G}_{1,C}) \otimes \ZZ_{(p)}\{\zeta\}
\end{gather*}
where the modules above have bases given by:
\begin{align*}
\mr{X}^\infty & := \bra{1_{sp^n/j,k}}, \quad p \not|s, n \ge 0, 1 \le k \le n+1, 1 \le j \le a_{n-k+1}, p^{k-1}|j,\\
& \quad \text{either} \: p^k \not| j \: \text{or} \: j > a_{n-k}, \\
\mr{X}_\infty^\infty & := \bra{1_{0/j,k}}, \quad  j \ge 1, k = \nu_p(j)+1, \\
\mr{X}\zeta^\infty_C & := \bra{\zeta_{sp^n/j,k}}, \quad p \not| s, n \ge 0: \\
& \quad \begin{cases}
\nu_p(s+1) = 0: & 1 \le k \le n+1, 1 \le j \le a_{n-k+1}, p^{k-1}|j, \\ 
& \text{either} \: p^k \not| j \: \text{or} \: j > a_{n-k}, \\
\nu_p(s+1) = i > 0: &    
\begin{cases}
1 \le k \le i-1: & 1 \le j \le a_{n-k+1}, p^{k-1}|j, \\ 
& \text{either} \: p^k \not| j \: \text{or} \: j > a_{n-k}, \\
i \le k \le n: & a_{n-k} < j \le a_{n-k+1}, p^k|j,
\end{cases}
\end{cases}\\
\mr{Y}^\infty_C & := \bra{(\mr{h}_1)_{sp/j,k}}, \quad 1 \le j < p-1, k= 1, \: \text{and} \: j = p-1, k = 2 \: \text{if} \: p|s, \\
\mr{Y}^\infty_{0,C} & := \bra{(\mr{h}_0)_{sp^n/j,k}}, \quad s \not\equiv 0, -1 \mod p, 1 \le k \le n, \\
& \quad A_{n-k}+2 < j \le A_{n-k+1}+2, p^{k-1}|j-1, \: \text{and} \: p^k|j-1 \: \text{if} \: j-1 \le a_{n-k+1}, \\ 
& \quad \text{as well as} \: j = 1, k = n+1. \\
\mr{Y}_{1,C}^\infty & := \bra{(\mr{h}_0)_{sp^n-p^{n-2}/j,k}}, \quad n \ge 2, s = p^ms', p \not| s', 1 \le k \le n+1: \\
& \begin{cases}
\quad p \not| j-1: & k = 1, a_{n-2}+1 < j \le p^n-p^{n-2}+A_{n-2}+2, \\ \\
\begin{array}{l}
p |j-1 \: \text{and} \\
j > p^{n}-p^{n-2}+1:
\end{array}  & 
\begin{array}{l}
k \ge 1, j = tp^{k-1}+1, \\
j \le p^n+p^{n-2}+A_{n-k-1}+2, \\ 
\text{and} \: p \not| t \: \text{or} \: j > p^n-p^{n-2}+A_{n-k-2}+2,
\end{array}\\ \\
\begin{array}{l}
p |j-1 \: \text{and} \\
j \le p^{n}-p^{n-2}+1:
\end{array} &
\begin{cases}
2 \le k \le n-2: & k \le m+1, \\ 
& j = tp^{k-1}+1, p \not| t, \\
& j > a_{n-k-1}+1, \\
k = n-1: & j = p^n-p^{n-2}+1, \: \\ 
& \text{or} \: j = 1 \: \text{and} \: n \le m+2, \\
k = n:  & j = tp^{n-1}-p^{n-2}+1, \\
&  n \le m+1, t \not\in \{p, p-1\}, \\
k = n+1: & j = p^n-p^{n-1}-p^{n-2}+1, \\ 
& n \le m,
\end{cases} 
\end{cases} \\
\mr{Y}^\infty_{\infty,C} & := \QQ/\ZZ_{(p)} \: \text{generated by} \: \{ \mr{h}_{0/1,k} \}, k \ge 1, \\
\mr{G}^\infty_C & := \bra{(\mr{G}_n)_{sp^n/j,k}}, \quad n \ge 0, 1 \le j \le a_n, s = p^{i}s', p \not| s' \\
& \quad \begin{cases} 
n = 0, s' \not\equiv -1 \mod p: & k = i+1, \\
n \ge 1, s' \not\equiv -1 \mod p: & k = \nu_p(j+A_{n-1}+1)+1 \le i+1, \\
n \ge 1, s' \equiv -1 \mod p: & k = \nu_p(j+A_{n-1}+1)+1 \le i,
\end{cases}\\
\mr{G}^\infty_0 & := \QQ/\ZZ_{(p)} \: \text{generated by} \: \{(\mr{G}_0)_{0/1,k}\}, k \ge 1, \\ 
\mr{Y}^{\infty,G}_{0, C} & := \bra{ (\mr{h}_0)_{sp^n/j,k}}, \quad n \ge 0, s \not\equiv 0, -1 \mod p, k \ge 1, j = tp^k+1, t \ne 0, \\
& \quad A_{n-k}+2 < j \le A_{n-k+1}+2, \\
\mr{Y}^{\infty, G}_{1,C} & := \bra{ (\mr{h}_0)_{sp^n-p^{n-2}/j,k})}, \quad n \ge 2, k \ge 1,  p^k|j+a_{n-1}, \\
& \quad p^n-p^{n-2}+A_{n-k-2}+2 < j \le p^n-p^{n-2}+A_{n-k-1}+2.
\end{align*}
\end{thm}

\begin{rmk}
Unlike in Theorem~\ref{thm:M20}, we have presented the modules in Theorem~\ref{thm:SY} in terms of an integral basis, as in \cite{SY}.  This way, the various modules are more easily compared to the corresponding modules in \cite{SY}. 
\end{rmk}

\begin{rmk}\label{rmk:error1}
The module $\mr{Y}^\infty_{1,C}$ differs from that which appears in Theorem~2.3 of \cite{SY} in two ways.  Firstly, the conditions ``$k \le m+1$'', ``$n \le m+2$'', ``$n \le m+1$'', and ``$n \le m$'' in the various subcases are absent from \cite{SY}.  These conditions are necessary, because they eliminate targets of differentials in the 
$v_0$-BSS (\ref{eq:v0BSS}).  The differentials in question are
$$  d(1)_{s'p^{n+m}/j+a_{n-1}} \doteq v_0^{m+1} (h_0)_{s'p^{n+m}-p^{n-2}/j} + \cdots $$
for $p \not| s'$, $j \le p^n-p^{n-2}$, $p^{m+1}| j+a_{n-1}$ (see Theorem~5.1 of \cite{SY}).  Secondly, in \cite{SY} the condition ``$j = tp^{k-1}+1$'' above instead reads ``$j = tp^{k}+1$''.  The source of this discrepancy is in Proposition~7.8 of \cite{SY}, where it is proven that there are differentials
$$ d((h_0)_{sp^n-p^{n-2}/j}) \doteq v_0^k(G_{n-k-1})_{sp^n-p^{n-1}/j-p^n+p^{n-2}-A_{n-k-2}-2} + \cdots $$
for $j \le p^{n}-p^{n-2}+A_{n-k-1}+2$ and $p^k|j+a_{n-1}$.  The issue is that the targets of these differentials are not present for $j \le p^{n}-p^{n-2}+A_{n-k-2}+2$.  While alternative targets are supplied by Proposition~7.8 of \cite{SY} for $j \le p^{n}-p^{n-2}+1$, the range $p^{n}-p^{n-2}+1 < j \le p^{n}-p^{n-2}+A_{n-k-2}+2$ is not addressed.  For the purposes of the projective $v_0$-BSS, however, Proposition~7.8 gives enough of a lower bound on the length of the projective $v_0$-BSS differential to deduce the orders of these groups in these missing cases.
\end{rmk}

\begin{rmk}\label{rmk:error2}
The module $\mr{G}^\infty_C$ differs from that which appears in Theorem~2.3 of \cite{SY} in three respects.  Firstly, in \cite{SY} there is the condition: 
\begin{quote}
``if $s' \not\equiv -1 \mod p$ then $p^{i+1} \not| j+A_{n-i-1}+1$. ''
\end{quote}
However, in light of Propositions~7.2 and 7.5 of \cite{SY}, this condition should instead read:
\begin{quote}
``if $s' \not\equiv -1 \mod p$ then $p^{i+1} \not| j+A_{n-1}+1$. ''
\end{quote}
Secondly, in \cite{SY} there is the condition:
\begin{quote}
``if $s' \equiv -1 \mod p^2$ then $p^{i} \not| j+A_{n-i}+1$. ''
\end{quote}
In light of Propositions~7.6 and 7.8 of \cite{SY}, this condition should instead read:
\begin{quote}
``if $s' \equiv -1 \mod p$ then $p^{i} \not| j+A_{n-1}+1$. ''
\end{quote}
Thirdly, the variable $i$ which appears in the second set of conditions describing $\mr{G}^\infty_C$ in Theorem 2.3 of \cite{SY} (i.e. the set of conditions involving the variable ``$l$'' in their notation) has nothing to do with the variable $i$ appearing in the first set of conditions describing $\mr{G}_C^\infty$.  This error arose because the definition of $\mr{G}^\infty_C$ at the top of page 287 of \cite{SY} involves superimposing the conditions of $\mr{G}_C$ on page 284 of \cite{SY}; both sets of conditions involve a variable ``$i$'', but these $i$'s are not the same.
\end{rmk}

\begin{rmk}\label{rmk:error3}
The module $\mr{Y}^{\infty,G}_{1,C}$ differs from that which appears in Theorem~2.3 of \cite{SY}.  We have replaced the condition
\begin{quote}
``$ p^k|j-1$''
\end{quote}
in \cite{SY} with the condition
\begin{quote}
``$ p^k|j+a_{n-1}$.''
\end{quote}
This only has the effect of adding the generators
$$ \mr{h}_0\zeta_{sp^n-p^{n-2}/p^n-p^{n-2}+1, n-1}. $$
These generators must be present, in light of Remark~9.10 of \cite{SY}, together with the $v_0$-BSS differential
$$ d (\mr{h}_0)_{sp^n-p^{n-2}/p^{n}-p^{n-2}+1} \doteq v_0^{n-1} (\mr{G}_0)_{sp^n-p^{n-1}/1}+ \cdots  $$
implied by Propositions~7.6 and 7.8 of \cite{SY}.
\end{rmk}

We give a dictionary between our presentation of $H^*M_0^2$ (Theorem~\ref{thm:M20}) and the Shimomura-Yabe presentation (Theorem~\ref{thm:SY}) below.  As before, our generators are italicized, while the Shimomura-Yabe generators are in non-italic typeface.  Family-by-family, we give a  \emph{basis} for our families, and then indicate the corresponding Shimomura-Yabe basis elements, broken down into cases.

\begin{align*}
X^\infty & = \mr{X}^\infty, \\
X_\infty^\infty & = \mr{X}^\infty_\infty, \\
Y_0^\infty & \ni (h_0)_{sp^n/j,k}, \quad s \not\equiv 0, -1 \, \mr{mod} \, p, n \ge 0, 1 \le k \le n+1, 2 \le j \le A_{n-k+1}+2, \\
& \quad p^{k-1}|j-1, \text{either} \: p^k \not| j-1 \: \text{or} \: j > A_{n-k}+2, \: \text{as well as} \: j = 1, k = n+1 \\
& = \begin{cases}
\zeta_{sp^n/j-1,k}, & 2 \le j \le a_{n-k+1}+1, \nu_p(j-1) = k-1,  \quad (\mr{X}\zeta^\infty_C) \\
(\mr{h}_0)_{sp^n/j,k}, & \text{either} \: a_{n-k+1} < j \le A_{n-k+1}+2, \nu_p(j-1) = k-1 \\
& \text{or} \: j > A_{n-k}+2 \: \text{or} \: j = 1, \quad(\mr{Y}^\infty_{0,C})
\end{cases}\\
Y_{0,\infty}^\infty  & \ni (h_0)_{0/j,k}, j \ge 2, k-1 = \nu_p(j-1) \: \text{and} \: \QQ/\ZZ_{(p)} \: \text{generated by} \: j = 1, k \ge 1, \\
& = 
\begin{cases}
\zeta_{0/j-1, k}, & j \ge 2, \quad (\mr{X}_\infty^\infty\{\zeta\}) \\
\mr{h}_{0/1,k}, & j = 1, \quad (\mr{Y}^\infty_{\infty, C})
\end{cases} \\
\zeta Y_{0,\infty}^\infty & = \mr{Y}^\infty_{\infty,C}\{\zeta\}, \\
Y^\infty & \ni  (h_1)_{sp/j,k}, \quad k = 1, 1 \le j < p-1, \: \text{and} \: j = p-1, k =  
\begin{cases}
1, & p \not| s, \\
2, & p| s
\end{cases} \\
& = 
\begin{cases}
(\mr{h}_1)_{sp/j,k}, & j < p-1 \: \text{and} \: j = p-1 \: \text{if} \: p|s, \quad (\mr{Y}^\infty_C) \\
\zeta_{sp/p, 1}, & j = p-1, p \not| s, \quad (\mr{X}\zeta^\infty_C)
\end{cases} \\
Y_1^\infty & \ni (h_0)_{sp^n-p^{n-2}/j,k}, \quad \text{writing} \: s = p^i s', p \not| s':\\
& \begin{cases}
j \le p^n-p^{n-2}:
& 
\begin{array}{l}
1 \le k \le \min(n+1, i+1), p^{k-1}|j+a_{n-1}, \\
\text{either} \: p^k \not| j+a_{n-1} \: \text{or}  \: k = i+1, 
\end{array}\\ \\
j > p^n-p^{n-2}:
& 
\begin{array}{l}
1 \le k \le n-1, j \le p^n-p^{n-2}+A_{n-k-1}+2, p^{k-1} | j+a_{n-1}, \\
\text{either} \: p^{k}\not| j+a_{n-1} \: \text{or} j > p^n-p^{n-2}+A_{n-k-2}+2 
\end{array}
\end{cases} \\
& =
\begin{cases}
\zeta_{sp^n/j+a_{n-1}, k}, & 1 \le j \le p^n-p^{n-2}, p^{k}|j+a_{n-1}, \quad (\mr{X}\zeta^\infty_C) \\
\zeta_{sp^n-p^{n-2}/j-1,k}, &  \nu_p(j+a_{n-1}) = k-1, j \le a_{n-k-1} + 1, \quad (\mr{X}\zeta^\infty_C)\\
(\mr{h}_0)_{sp^n-p^{n-2}/j,k}, & \text{otherwise}, \quad (\mr{Y}^\infty_{1,C})
\end{cases} \\
G^\infty & \ni (G_n)_{sp^n/j,k}, \quad n \ge 0, 1 \le j \le a_n, \: \text{writing} \: s = p^i t, p \not| t, \: \text{we have}: \\
& \quad \begin{cases}
t \not\equiv -1 \mod p: & i\ge 0, 
\begin{cases}
n = 0: & k = i+1, \\
n \ge 1: &  k = \min(\nu_p(j+A_{n-1}+1)+1, i+1),
\end{cases}
\\
t \equiv -1 \mod p: & i \ge 1,
\begin{cases}
n = 0: & k = i, \\
n \ge 1:    
& k = \min(\nu_p(j+A_{n-1}+1)+1, i),
\end{cases}
\end{cases} \\
& =
\begin{cases}
(\mr{G}_0)_{s/1,i+1}, & n = 0, t \not\equiv -1 \mod p, \quad (\rm{G}^\infty_C) \\
\mr{h}_0\zeta_{t'p^{i+1}-p^{i-1}/p^{i+1}-p^{i-1}+1, i},  & n = 0, t = t'p-1, \quad (\mr{Y}^{\infty,G}_{1,C}\{\zeta\})  \\
(\mr{G}_n)_{sp^n/j,k}, & n \ge 1, p^k \not| j+A_{n-1}+1, \quad (\mr{G}^\infty_C) \\
\mr{h}_0\zeta_{tp^{n+i}/j+A_{n-1}+2, k}, & n \ge 1, t \not\equiv -1 \mod p, p^k|j+A_{n-1}+1, \\
& \hfill (\mr{Y}^{\infty,G}_{0,C}\{\zeta\}) \\
\mr{h}_0\zeta_{\frac{t'p^{n+i+1}-p^{n+i-1}}{j+p^{n+i+1}-p^{n+i-1}+A_{n-1}+2,k}}, & n \ge 1, t = t'p-1, p^k|j+A_{n-1}+1, \\
& \hfill (\mr{Y}^{\infty,G}_{1,C} \{\zeta\} )
\end{cases} \\
G_\infty^\infty & \ni (G_n)_{0/j,k}, \quad n \ge 0, 1 \le j \le a_n, 
\begin{cases} 
n = 0: & \text{generates} \: \QQ/\ZZ_{(p)}, k \ge 1, \\
n > 0: & 1 \le k \le n+1, 1 \le j \le a_n, \\
& k = \nu_p(j+A_{n-1}+1)+1
\end{cases}\\
& = 
\begin{cases}
(\mr{G}_0)_{0/1,k},  & n = 0, \quad (\mr{G}_0^\infty) \\
(\mr{G}_n)_{0/j,k}, & n \ge 1, \quad (\mr{G}^\infty_C) \\
\end{cases}\\
\zeta G_\infty^\infty & = \mr{G}^\infty_0 \{\zeta\}. \\
\end{align*}

\section{$E(2)$ and $K(2)$-local computations}\label{sec:K(2)S}

The computation of the groups $\pi_*M(p)_{E(2)}$,  $\pi_*M(p)_{K(2)}$, $\pi_* S_{E(2)}$ and $\pi_* S_{K(2)}$ follow quickly from $H^*M_1^1$ and $H^*M_0^2$.  We briefly review this in this section.

The Morava change of rings theorem, applied in the context of $n = 0$, gives the following well known fact.

\begin{lem}
We have 
$$
H^{s,t}M_0^0 \cong \begin{cases}
\QQ, & (s,t) = (0,0), \\
0, & \text{otherwise}.
\end{cases}
$$
\end{lem}

\begin{thm}[Theorem~1.2 of \cite{Ravenel}]
We have
$$ H^{s,t}M_1^0 \cong \FF_p[v_1^{\pm 1}] \otimes E[h_0] $$
where 
\begin{align*}
\abs{v_1} & = (0,q), \\
\abs{h_0} & = (1,q).
\end{align*}
\end{thm}

In the following theorem, we are using the notation
$$ x_{s/k} := p^{-k} v_1^s x, \quad \text{for} \: x \in H^*(M_1^0) $$
to refer to elements in $H^*M_0^1$.

\begin{thm}[Theorem~4.2 of \cite{MillerRavenelWilson}]
The groups $H^*M_0^1$ are spanned by
\begin{align*}
1_{s/k}, & \quad k \ge 1, p^{k-1}|s, \\
(h_0)_{-1/k}, & \quad k \ge 1 \\
\end{align*}
\end{thm}

The ANSS's
\begin{align*}
H^{s,t}M_0^0 & \Rightarrow \pi_{t-s} M_0(S) \\
H^{s,t}M_1^0 & \Rightarrow \pi_{t-s} M_1(M(p)) \\
H^{s,t}M_0^1 & \Rightarrow \pi_{t-s-1} M_1(S) \\
H^{s,t}M_1^1 & \Rightarrow \pi_{t-s-1} M_2(M(p)) \\
H^{s,t}M_0^2 & \Rightarrow \pi_{t-s-2} M_2(S) 
\end{align*}
all collapse because of their sparsity.  

Consider the chromatic spectral sequence
$$ E_1^{n,k} = \bigoplus_{n = 1}^2 \pi_k M_n(M(p)) \Rightarrow \pi_k M(p)_{E(2)}. $$
The differentials are given by
\begin{align*}
d_1(1_{s}) & = 
\begin{cases}
1_{0/-s}, & s < 0, \\
0, & s \ge 0,
\end{cases} \\
d_1((h_0)_s) & = 
\begin{cases}
(h_0)_{0/-s}, & s < 0, \\
0, & s \ge 0.
\end{cases}
\end{align*}
We therefore get the following well-known consequence of Shimomura's calculation of $H^*M_1^1$.  Here, the degrees of the elements are their internal degrees, viewed as elements of $H^*M_i^j$, and the homological grading is to be ignored.

\begin{thm}
We have
\begin{multline*}
\pi_* M(p)_{E(2)} \cong \FF_p[v_1]\otimes E[h_0] \oplus 
(\Sigma^{-1}X_{\infty} \oplus  \Sigma^{-2}Y_\infty)\{\zeta\} \oplus \\
(\Sigma^{-1} X \oplus \Sigma^{-2} (Y_0 \oplus Y \oplus Y_1)
 \oplus \Sigma^{-3} G) \otimes E[\zeta] 
\end{multline*}
where $\abs{\zeta} = -1$.
\end{thm}

Using the $\lim^i$ sequence associated to 
$$ M(p)_{K(2)} \simeq \holim_j M(p, v_1^j)_{E(2)} $$
we get the following theorem (see Section~15.2 of \cite{HoveyStrickland}).

\begin{thm}\label{thm:M(p)K(2)}
We have
\begin{multline*}
\pi_* M(p)_{K(2)} \cong \FF_p[v_1]\otimes E[h_0, \zeta]  \oplus
(\Sigma^{-1} X \oplus \Sigma^{-2} (Y_0 \oplus Y \oplus Y_1)
 \oplus \Sigma^{-3} G) \otimes E[\zeta] 
\end{multline*}
where $\abs{\zeta} = -1$.
\end{thm}

Consider the chromatic spectral sequence
$$ E_1^{n,k} = \bigoplus_{n = 0}^2 \pi_k M_n(S) \Rightarrow \pi_k S_{E(2)}. $$
The differential
$$ d_1: \QQ = \pi_0 M_0(S) \rightarrow \pi_{-1} M_1(S) = \QQ/\ZZ_{(p)} \bra{(h_0)_{-1/k}\: : \: k \ge 1} $$
is the canonical surjection.
The differentials
$$ d_1: \pi_k M_1(S) \rightarrow \pi_{k-1} M_2(S) $$
are given by
\begin{align*}
d_1(1_{s/k}) & = 
\begin{cases}
1_{0/-s,k}, & s < 0, \\
0, & s \ge 0,
\end{cases}\\
d_1((h_0)_{-1/k}) & = (h_0)_{0/1,k}.
\end{align*}

Write
\begin{align*}
Y^\infty_{0,\infty} & = Y^\infty_{0,\infty}[0] \oplus Y^\infty_{0,\infty}[1], \\
G_\infty^\infty & = G_{\infty}^\infty[0] \oplus G^\infty_{\infty}[1]
\end{align*}
where
\begin{align*}
Y_{0,\infty}^\infty[0] & = \bra{(h_0)_{0/1,k}\: : \: k \ge 1}, \\
Y^\infty_{0,\infty}[1] & = \bra{(h_0)_{0/j,k}\: : \: j \ge 2, p^{k-1}|j-1}, \\
G_\infty^\infty[0] & = \bra{(G_0)_{0/1,k} \: : \: k \ge 1}, \\
G_\infty^\infty[1] & = \bra{(G_n)_{0/j,k} \: : \: n \ge 1, 1 \le j \le a_n, p^{k-1}|j+A_{n-1}+1}.
\end{align*}
We deduce the following main theorem of \cite{SY}.  

\begin{thm}[Theorem~2.4 of \cite{SY}]
We have 
\begin{align*}
\pi_* S_{E(2)} \cong \: & \ZZ_{(p)} \oplus \Sigma^{-1} \bra{1_{sp^n/n+1} \: : \: n \ge 0, s > 0, p \not| s} \\
&  \Sigma^{-2} X^{\infty} \oplus \Sigma^{-3} (Y_0^\infty \oplus Y^\infty_{0,\infty}[1] \oplus Y^\infty \oplus Y_1^\infty) \oplus \\
& \Sigma^{-4}(\zeta Y^\infty_{0,\infty} \oplus G^\infty \oplus G^\infty_\infty) \oplus \Sigma^{-5} \zeta G_\infty^\infty.
\end{align*}
\end{thm}

Using the $\lim^i$ sequence associated to 
$$ S_{K(2)} \simeq \holim_{j,k} M(p^k, v_1^j)_{E(2)} $$
we get the following theorem.

\begin{thm}\label{thm:K(2)S}
We have 
\begin{align*}
\pi_* S_{K(2)} \cong \: & \ZZ_{p}\otimes E[\zeta, \rho] \oplus \Sigma^{-1} \bra{1_{sp^n/n+1} \: : \: n \ge 0, s > 0, p \not| s} \otimes E[\zeta] \oplus \\
&  \Sigma^{-2} X^{\infty} \oplus \Sigma^{-3} (Y_0^\infty \oplus Y^\infty \oplus Y_1^\infty) \oplus \Sigma^{-4}(G^\infty \oplus G^\infty_\infty[1])
\end{align*}
where $\abs{\zeta} = -1$ and $\abs{\rho} = -3$. 
\end{thm}

\begin{rmk}
The existence of the exterior algebra factors involving $\zeta$ and $\rho$ in Theorem~\ref{thm:K(2)S} are closely related to Hopkins' chromatic splitting conjecture (see \cite{Hovey}).  In fact, using the fiber sequence
$$ M_2(S) \rightarrow S_{K(2)} \rightarrow S_{K(2),E(1)} $$
one easily deduces
\begin{multline*} \pi_* S_{K(2),E(1)} \cong \\ 
(\ZZ_p \oplus \Sigma^{-1}\bra{1_{sp^n/n+1} \: : \: n \ge 0, p \not|s} \oplus \Sigma^{-2} \QQ/\ZZ_{(p)}) \otimes E[\zeta] \oplus \Sigma^{-3} \QQ_p \oplus \Sigma^{-4} \QQ_p,
\end{multline*}
as predicted by the chromatic splitting conjecture.
\end{rmk}

\section{Gross-Hopkins duality}\label{sec:GHD}

The reader may notice that the patterns which occur in Figure~\ref{fig:M11} are ambigrammic: they are invariant under rotation by $180^\circ$.  This is explained by Gross-Hopkins duality.

To proceed, we must work with Picard group graded homotopy.
The following is an unpublished result of Hopkins.

\begin{thm}[Hopkins]
There is an isomorphism 
\begin{equation}\label{eq:Pic}
\mr{Pic}_{K(2)} \cong \ZZ_p \times \ZZ_p \times \ZZ/2(p^2-1).
\end{equation}
The group is topologically generated by $S^1_{K(2)}$ and $S^0_{K(2)}[\det]$.  The isomorphism (\ref{eq:Pic}) can be chosen so that these generators are given by
\begin{align}
S^1_{K(2)} & = (1, 0, 1), \label{eq:S1} \\
S^0_{K(2)}[\det] & = (0,1,2(p+1)). \label{eq:Sdet}
\end{align}
\end{thm}

\begin{proof}[Overview of the proof]
As this isomorphism is not in print, we give a brief explanation (note that the analogous fact for $p = 3$ is published, see \cite{Karamanov}).  Given an object $X \in \mr{Pic}_{K(2)}$, the associated Morava module $(E_2)^\wedge_* X$ is invertible.  In particular, as a graded $(E_2)_*$-module, it is free of rank $1$, concentrated either in even or odd degrees.  Define $\epsilon(X)  \in \ZZ/2$ to be the degree of a generator of $(E_2)^\wedge_*X$.  This gives a short exact sequence
\begin{equation}\label{eq:parity}
0 \rightarrow \mr{Pic}_{K(2)}^0 \xrightarrow{\iota_0} \mr{Pic}_{K(2)} \xrightarrow{\epsilon} \ZZ/2 \rightarrow 0.
\end{equation}
Since invertible Morava modules are in bijective correspondence with degree $1$ group cohomology classes, taking the degree zero part of the associated Morava module gives a map
\begin{equation}\label{eq:map1}
 \mr{Pic}_{K(2)}^0 \xrightarrow{(E_2)^\wedge_0(-)} H^1_c(\GG_2; (E_2)_0^\times) \cong H^1_c(\MS_2; (E_2)_0^\times)^{\mit{Gal}}.
\end{equation}
(Here, $\mit{Gal}$ denotes the Galois group of $\FF_{p^2}/\FF_p$.)  Since the reduction map
$$ (E_2)_0 \cong \WW[[u_1]] \rightarrow \WW $$
is equivariant with respect to the subgroup $\WW^\times < \MS_2$ (where $\WW$ denotes the Witt ring of $\FF_{p^2}$), there is a map
\begin{equation}\label{eq:map2}
H^1_c(\MS_2; (E_2)_0^\times)^{\mit{Gal}} \rightarrow H_c^1(\WW^\times; \WW^\times)^{\mit{Gal}} \cong \mr{End}^c(\WW^\times)^{\mit{Gal}}.
\end{equation}
The crux of Hopkins' argument is that both (\ref{eq:map1}) and (\ref{eq:map2}) are isomorphisms, and there is an isomorphism
$$ \mr{End}^c(\WW^\times)^{\mit{Gal}} \underset{(\dag)}{\cong} \ZZ_p \times \ZZ_p \times \ZZ/(p^2-1). $$
The isomorphism $(\dag)$ follows from the usual Galois-equivariant isomorphism
$$ \WW \times \FF_{p^2}^\times \xrightarrow[\cong]{\exp(px) \times \tau} \WW^\times $$
where $\tau$ is the Teichm\"uller lift.
Since there are no continuous group homomorphisms between $\FF^\times_{p^2}$ and $\WW$, we get
$$ \mr{End}^c(\WW^\times)^{\mit{Gal}} \xrightarrow{\cong} \mr{End}^c(\WW)^{\mit{Gal}} \times \mr{End}(\FF_{p^2}^\times)^{\mit{Gal}}. $$
Every endomorphism of $\FF_{p^2}^\times$ is Galois equivariant (since the Galois action is the $p$th power map), and we have
$$ \mr{End}(\FF_{p^2}^\times) \cong \ZZ/(p^2-1). $$
There is an isomorphism
$$ \mr{End}^c(\WW)^{\mit{Gal}} \cong \ZZ_p\{\mr{Id}, \mr{Tr}\}. $$
The Galois equivariant endomorphism of $\WW^\times$ induced from $[S^2_{K(2)}] \in \mr{Pic}^0_{K(2)}$ (respectively $[S^0_{K(2)}[det]] \in \mr{Pic}^0_{K(2)}$) is the identity (respectively the norm).  It follows that under isomorphisms (\ref{eq:map1}), (\ref{eq:map2}), and $(\dag)$ above, we have:
\begin{align*}
S^2_{K(2)} & = (1, 0, 1), \\
S^0_{K(2)}[\det] & = (0,1,p+1).
\end{align*}
Since $\epsilon[S^1_{K(2)}] = 1$ and $2[S^1_{K(2)}] = [S^2_{K(2)}]$ in $\mr{Pic}_{K(2)}$, we deduce from (\ref{eq:parity}) isomorphism (\ref{eq:Pic}).  Moreover, the induced map
$$ \ZZ_p \times \ZZ_p \times \ZZ/(p^2-1) \cong \mr{Pic}^0_{K(2)} \hookrightarrow \mr{Pic}_{K(2)} \cong \ZZ_p \times \ZZ_p \times \ZZ/2(p^2-1) $$
can be taken to be $(a, b, c) \mapsto (2a, b, 2c)$.  The identities (\ref{eq:S1}) and (\ref{eq:Sdet}) follow.
\end{proof}

The isomorphism (\ref{eq:Pic}) implies that we can $K(2)$-locally $p$-adically interpolate the spheres to get
\begin{gather}
S_{K(2)}^{s\abs{v_2}+i} = (s\abs{v_2}+i, 0, i),  \quad \text{for} \: s \in \ZZ_p, 0 \le i < 2(p^2-1),  \\
S_{K(2)}^{(1+p+p^2+ \cdots)\abs{v_2}+q+4} = (0, 0, 2(p+1)). \label{eq:2(p+1)} 
\end{gather}
For a $K(2)$-local spectrum $X$, we may define $\pi_{*,*}(X)$ by
$$ \pi_{s\abs{v_2}+i, j}(X) := [S_{K(2)}^{s\abs{v_2}+i}[{\det}^{j}], X] $$
for $s,j \in \ZZ_p$, $0 \le i < \abs{v_2}$.

By extending the families described in Theorems~\ref{thm:M11} and \ref{thm:M20} to allow for $s$ to lie in $\ZZ_p$ instead of $\ZZ$, one can regard Theorems~\ref{thm:M(p)K(2)} and \ref{thm:K(2)S} as giving $\pi_{*,0} M(p)_{K(2)}$ and $\pi_{*,0} S_{K(2)}$, where $*$ varies $p$-adically.  The author does not know how to compute $\pi_{*,j} S_{K(2)}$ for arbitrary $j \in \ZZ_p$.  However, as the following proposition illustrates, after smashing with the Moore spectrum $M(p)$ the elements $(a,*,b) \in \mr{Pic}_{K(2)}$ (under the isomorphism (\ref{eq:Pic})) are all equivalent for fixed $a$ and $b$ and $*$ ranging through $\ZZ_p$.

\begin{prop}
\begin{equation}\label{eq:M(p)[det]}
M(p)_{K(2)}[\det] \simeq \Sigma^{(1+p+p^2+ \cdots)\abs{v_2}+q+4}M(p)_{K(2)}.
\end{equation}
\end{prop}

\begin{proof}
Since the mod $p$ determinant takes values in $\FF_p^\times$, there is an isomorphism of Morava modules
$$ (E_2)^\wedge_* M(p) [{\det}^{p-1}] \cong (E_2)^\wedge_* M(p). $$
It follows that under isomorphism (\ref{eq:Pic}), the subgroup of $\ZZ_p \times \ZZ_p \times \ZZ/2(p^2-1)$ generated by $(0,p-1,0)$ acts trivially on $M(p)_{K(2)}$.  Thus the element in $\mr{Pic}_{K(2)}$ corresponding to $(0,1,0)$ also acts trivially.  The proposition follows from (\ref{eq:Sdet}) and (\ref{eq:2(p+1)}).
\end{proof}

Following \cite{GrossHopkins}, we define
$$ I_2X := IM_2(X) $$
where $I$ denotes the Brown-Comenetz dual.
The following proposition explains the self-duality apparent in Figure~\ref{fig:M11}.

\begin{prop}\label{prop:duality}
There is an equivalence
$$ I_2 M(p) \simeq \Sigma^{(1+p+p^2+\cdots)\abs{v_2} + q + 5}M(p)_{K(2)}. $$
\end{prop}


\begin{proof}
Theorem~6 of \cite{GrossHopkins}, when specialized to our case, states that there is an equivalence:
\begin{equation}\label{eq:GHD}
I_2 S \simeq S_{K(2)}^2[\det].
\end{equation}
Smashing (\ref{eq:GHD}) with $M(p)$ and using (\ref{eq:M(p)[det]}) we get
\begin{align*}
I_2 M(p) & \simeq \Sigma^{-1} M(p) \wedge I_2 S \\
& \simeq \Sigma^{-1} M(p) \wedge S_{K(2)}^2[\det] \\
& \simeq \Sigma^{(1+p+p^2+\cdots)\abs{v_2} + q + 5}M(p)_{K(2)}.
\end{align*}
\end{proof}

Unfortunately, as we have not given a method to compute $\pi_{*,j} S_{K(2)}$ for arbitrary $j$, (\ref{eq:GHD}) gives little insight into the shifted self-duality present in the patterns shown in Figure~\ref{fig:M20}.  However, using (\ref{eq:GHD}), one can turn the patterns of Figure~\ref{fig:M20} $180^\circ$ and regard them as being descriptions of the corresponding patterns occurring in the homotopy of $S^0_{K(2)}[\det]$.

\begin{rmk}
One way to compute the portion of $\pi_{*,j} S_{K(2)}$ spanned by elements of Adams-Novikov filtration $2$ is to adapt the method of congruences of modular forms of \cite{Behrens} to the situation: one just needs to twist the operators acting on the modular forms by appropriate powers of the determinants of the corresponding elements of $GL_2(\QQ_\ell)$.  In fact, this method helped the author correct an additional family of errors in $Y^\infty_1$ and $G^\infty$ which he missed in an earlier version of this paper. 
\end{rmk}

\section{A simplified presentation}\label{sec:M20new}

The patterns of Figure~\ref{fig:M20} suggest that we may reorganize the families $X$, $Y$, $Y_0$, $Y_1$, $G$, into four simple families, as explained in the following theorem.  In the theorem below, we have
$$ \abs{x(j,k)_{s}} = \abs{x} + s\abs{v_2} - jq. $$
We warn that while such an element $x(j,k)_s$ does have order $p^k$, the $j$ in the notation is not intended to indicate anything about $v_1$-multiplication.

\begin{thm}\label{thm:M20new}
$H^*M_0^2$ admits the following alternate presentation.
$$ H^*M_0^2 \cong X^\infty \oplus Y(0)^\infty \oplus Y(1)^\infty \oplus G^\infty \oplus X_\infty^\infty \oplus Y(0)^\infty_\infty \oplus \zeta Y(0)^\infty_\infty \oplus G_\infty^\infty \oplus \zeta G_\infty^\infty $$
where
\begin{align*}
X^\infty & := \bra{1(j,k)_{sp^n}}, \quad p \not|s, n \ge 0, 1 \le k \le n+1, 1 \le j \le a_{n-k+1}, p^{k-1}|j, \\
Y(0)^\infty & := \bra{h_0(j,k)_{sp^n}}, \quad p \not| s, \\
& \quad \begin{cases}
s \not\equiv -1 \mod p: & n \ge 0, 1 \le k \le n+1, 1 \le j \le A_{n-k+1}+2, \\ 
& p^{k-1}|j-1, \\
s \equiv -1 \mod p: & n \ge 1, 1 \le k \le n, 1 \le j \le A_{n-k}+2, \\ 
& p^{k-1} | j-1,
\end{cases}\\
Y(1)^\infty & := \bra{h_1(j,k)_{sp^n}}, \quad p \not| s, n \ge 1, 1 \le k \le n, 2 \le j+1 \le a_{n-k+1}, p^{k-1}|j+1, \\
G^\infty & := \bra{G_i(j,k)_{sp^n}}, \quad p\not| s, \\
& \quad \begin{cases}
s \not\equiv -1 \mod p: & n \ge 0, 0 \le i \le n, 1 \le j \le a_i, \\
& 1 \le k \le \min (i+1, n-i+1), p^{k-1}|j+A_{i-1}+1, \\
& (1 \le k \le n+1 \: \text{if} \: i = 0), \\
s \equiv -1 \mod p: & n \ge 1, 0 \le i \le n-1, 1 \le j \le a_i, \\
& 1 \le k \le \min (i+1, n-i), \: p^{k-1}|j+A_{i-1}+1, \\
& (1 \le k \le n \: \text{if} \: i = 0),
\end{cases} \\
X_\infty^\infty & := \bra{1(j,k)_0}, \quad k \ge 1, j \ge 1, p^{k-1}|j, \\
Y(0)^\infty_\infty & := \bra{h_0(j,k)_0}, \quad k \ge 1, j \ge 1, p^{k-1}|j-1, \\
\zeta Y(0)^\infty_\infty & := \bra{\zeta h_0(1,k)_0}, \quad k \ge 1, \\
G_\infty^\infty & := \bra{G_i(j,k)_0}, \quad i \ge 0, 1 \le j \le a_i, 1 \le k \le i+1, p^{k-1}|j+A_{i-1}+1, \\
& \quad (1 \le k \le \infty \: \text{if} \: i = 0), \\
\zeta G_\infty^\infty & := \bra{\zeta G_0(1,k)_0}, \quad k \ge 1.
\end{align*}
\end{thm}

\begin{figure}
\includegraphics[height=\textwidth, angle=270]{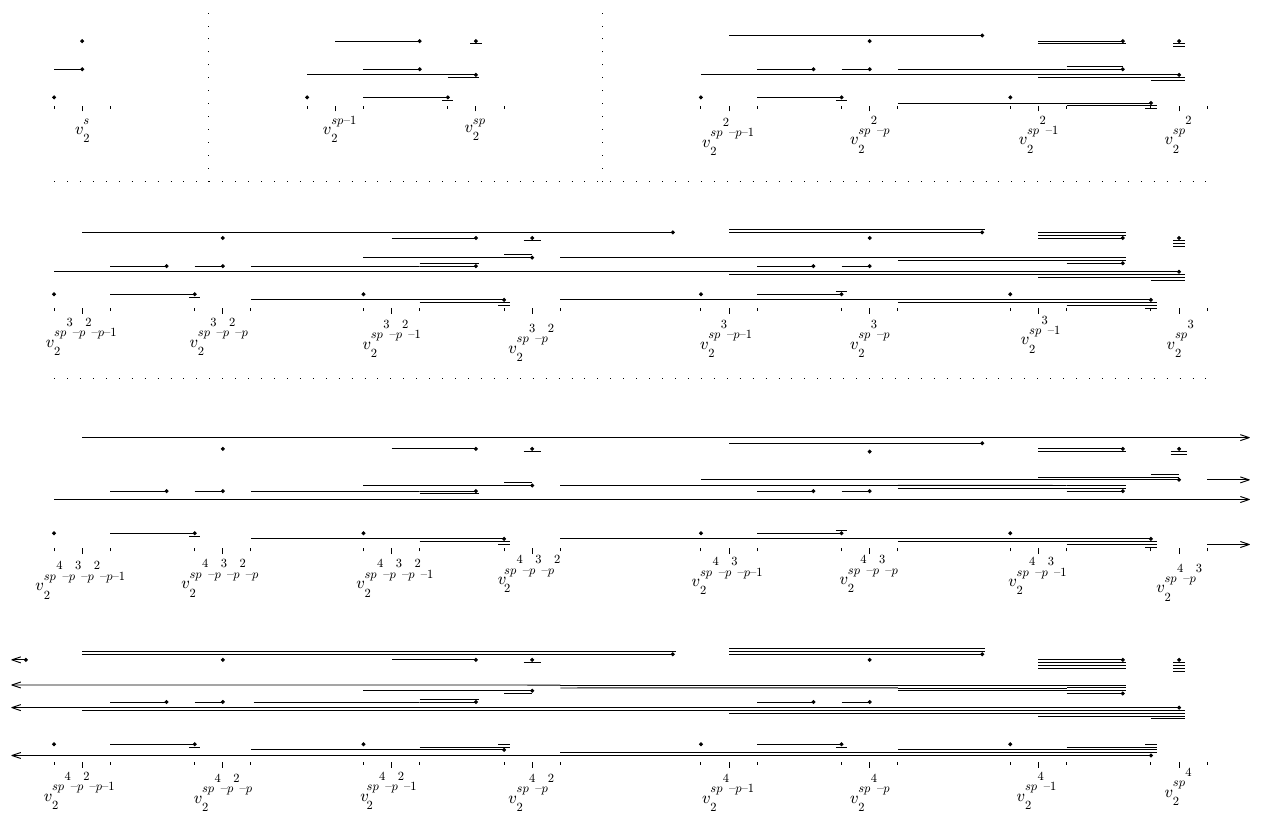}
\caption{$H^* M_0^2$ in the vicinity of $v_2^{sp^n}$, $0 \le n \le 4$, $s \not\equiv 0, -1 \mod p$ with respect to the simplified presentation}\label{fig:M20new}
\end{figure}

Figure~\ref{fig:M20new} shows the resulting patterns in the vicinities of $v_2^{sp^n}$ for $s \not\equiv -1 \mod p$ and $n \le 4$.  The meaning of the notation is identical to that of Figure~\ref{fig:M20} except that the lines are serving as an organizational principle, and are no longer meant to necessarily imply $v_1$-multiplication.

In order to prove that the presentation of Theorem~\ref{thm:M20new} is valid, we must provide a dictionary between the presentation of Theorem~\ref{thm:M20new} and the presentation of Theorem~\ref{thm:M20}.  
The modules
$$ X^\infty, X^\infty_\infty, G^\infty, G^\infty_\infty, \zeta G^\infty_\infty $$
share the same notation and indeed refer to the same modules as in Theorem~\ref{thm:M20}, with
$$ x(j,k)_s = x_{s/j,k}. $$
We also have 
\begin{align*}
Y_{0,\infty}^\infty & = Y(0)^\infty_\infty, \\
\zeta Y^\infty_{0, \infty} & = \zeta Y(0)^\infty_\infty.
\end{align*}
However, the modules $Y_0^\infty$, $Y^\infty$, and $Y_1^\infty$ of Theorem~\ref{thm:M20} get reorganized into the modules $Y(0)^\infty$ and $Y(1)^\infty$ of Theorem~\ref{thm:M20new}:
\begin{align*}
Y(0)^\infty \ni h_0(j,k)_{sp^n} & = 
\begin{cases}
(h_0)_{sp^n/j,k},  & s \not \equiv -1 \mod p,  \quad (Y_0^\infty) \\
(h_0)_{sp^n+p^n-p^{n-1}/j+p^{n+1}-p^{n-1}, k}, & s \equiv -1 \mod p, \quad (Y_1^\infty) \\
\end{cases} \\
Y(1)^\infty \ni h_1(j,k)_{sp^n} & = 
\begin{cases}
(h_1)_{sp^n/j,k},  & a_0 < j+1 \le a_1, \quad (Y^\infty) \\
(h_0)_{sp^n-p^{n-i}/j-a_{i-1}+1,k} & a_{i-1} < j+1 \le a_i, i > 1. \quad (Y_1^\infty)
\end{cases}
\end{align*}

The advantage of the presentation of Theorem~\ref{thm:M20new} is that it attaches to \emph{every} element $v_2^{sp^n}$
four $v_1$-torsion families: the two ``unbroken'' families $X^\infty$ and $Y(0)^\infty$ and the two ``broken'' families $Y(1)^\infty$ and $G^\infty$.  The unbroken families behave uniformly in $s$ and $n$, whereas the broken families display an exceptional behavior when $s \equiv -1 \mod p$.  This allows for easy understanding of the structure of $H^{s,t}M_0^2$ for $t \le 0$. The torsion bounds on $X^\infty$ and $Y(1)^\infty$ match up, as do the torsion bounds on $Y(0)^\infty$ and $G^\infty$.    Moreover, each of the four families are no more complicated than $X^\infty$, which corresponds to the family $\beta_{i/j,k}$ of \cite{MillerRavenelWilson}.  In contrast the presentation of $Y_1^\infty$ in Theorem~\ref{thm:M20} has a more complex feel to it, and the presentation of $\mr{Y}_1^\infty$ in Theorem~\ref{thm:SY} borders on incomprehensible.

The disadvantages of the presentation of Theorem~\ref{thm:M20new} is that we have forsaken a complete description of $v_1$-multiplication between the generators.  We have also broken any semblance of the Gross-Hopkins self-duality that was so readily apparent in Figure~\ref{fig:M11}.


\end{document}